%%%%%%%%%%%%%%%%%%%%%%%%%%%%%%%%%%%%%%%%%%%%%%%%%%%%%%%%%%
%%%%%%%%%%%%%%%%%%%%%%%%%%%%%%%%%%%%%%%%%%%%%%%%%%%%%%%%%%
%%
%%     This is the AMS-LaTeX file:
%%
%%     Colli-Gilardi-Sprekels 19
%%     Optimal distributed control of a generalized
%%     fractional Cahn--Hilliard system
%%     
%%%%%%%%%%%%%%%%%%%%%%%%%%%%%%%%%%%%%%%%%%%%%%%%%%%%%%%%%%

%\def\TeX{\input /articoli/ltx-tex/macrotex }

\def\LaTeX{%
  \let\Begin\begin
  \let\End\end
  \let\salta\relax
  \let\finqui\relax
  \let\futuro\relax}

\def\UK{\def\our{our}\let\sz s}
\def\USA{\def\our{or}\let\sz z}

\UK
%\USA

%%%%%%%%%%%%%%%%%%%%%%%%%%%%%%%%%

% scegliere fra \TeX e \LaTeX  e fra  \UK oppure \USA

%\TeX
\LaTeX

%\UK
\USA

%%%%%%%%%%%%%%%%%%%%%%%%%%%%%%%%%
%% page layout
%%%%%%%%%%%%%%%%%%%%%%%%%%%%%%%%%

\salta

\documentclass[twoside,12pt]{article}
\setlength{\textheight}{24cm}
\setlength{\textwidth}{16cm}
\setlength{\oddsidemargin}{2mm}
\setlength{\evensidemargin}{2mm}
\setlength{\topmargin}{-15mm}
\parskip2mm

%%%%%%%%%%%%%%%%%%%%%%%%%%%%%%%%%
%% packages
%%%%%%%%%%%%%%%%%%%%%%%%%%%%%%%%%

%\usepackage{color}
\usepackage[usenames,dvipsnames]{color}
\usepackage{amsmath}
\usepackage{amsthm}
\usepackage{amssymb}
\usepackage[mathcal]{euscript}
\usepackage{cite}

%\usepackage[notref,notcite]{showkeys}
%\usepackage{showkeys}
%
%		COLORS FOR CORRECTIONS
%
% do the same, please (i.e., don't use the standard {\color{red} text} or similar): 
% just choose the color you prefer in \def\yourname

% example of use:  \juerg{I want this to become blue}

\definecolor{viola}{rgb}{0.3,0,0.7}
\definecolor{ciclamino}{rgb}{0.5,0,0.5}

\def\pier #1{{\color{red}#1}}
\def\giapie #1{{\color{blue}#1}}
\def\juerg #1{{\color{Green}#1}}

\def\revis #1{{\color{red}#1}}
\def\revis #1{#1}

\def\pier #1{#1}
\def\juerg #1{#1}
\def\oldgianni #1{#1}
\def\giapie #1{#1}

%%%%%%%%%%%%%%%%%%%%%%%%%%%%%%%%%
%% you may adjust the baseline
%%%%%%%%%%%%%%%%%%%%%%%%%%%%%%%%%

%\renewcommand{\baselinestretch}{0.95}

%%%%%%%%%%%%%%%%%%%%%%%%%%%%%%%%%
%% bibliographystyle
%%%%%%%%%%%%%%%%%%%%%%%%%%%%%%%%%

\bibliographystyle{plain}

%%%%%%%%%%%%%%%%%%%%%%%%%%%%%%%%%
%% environments
%%%%%%%%%%%%%%%%%%%%%%%%%%%%%%%%%

%

\finqui

\def\Beq{\Begin{equation}}
\def\Eeq{\End{equation}}
\def\Bsist{\Begin{eqnarray}}
\def\Esist{\End{eqnarray}}

\def\Bthm{\Begin{theorem}}
\def\Ethm{\End{theorem}}
\def\Blem{\Begin{lemma}}
\def\Elem{\End{lemma}}
\def\Bprop{\Begin{proposition}}
\def\Eprop{\End{proposition}}
\def\Bcor{\Begin{corollary}}
\def\Ecor{\End{corollary}}
\def\Brem{\Begin{remark}\rm}
\def\Erem{\End{remark}}

\def\Bdim{\Begin{proof}}
\def\Edim{\End{proof}}
\def\Bcenter{\Begin{center}}
\def\Ecenter{\End{center}}
\let\non\nonumber

%%%%%%%%%%%%%%%%%%%%%%%%%%%%%%%%%
%% macros
%%%%%%%%%%%%%%%%%%%%%%%%%%%%%%%%%

% macro salvate

% sottosezioni non numerate

\def\step #1 \par{\medskip\noindent{\bf #1.}\quad}

% abbreviazioni di parole

\def\aand{\quad\hbox{and}\quad}

\def\lhs{left-hand side}
\def\rhs{right-hand side}
\def\sfw{straightforward}

% versioni inglesi (UK) o americane (USA)

% bold, cal e mathop

\def\multibold #1{\def\arg{#1}%
  \ifx\arg\pto \let\next\relax
  \else
  \def\next{\expandafter
    \def\csname #1#1#1\endcsname{{\bf #1}}%
    \multibold}%
  \fi \next}

\def\pto{.}

\def\multical #1{\def\arg{#1}%
  \ifx\arg\pto \let\next\relax
  \else
  \def\next{\expandafter
    \def\csname cal#1\endcsname{{\cal #1}}%
    \multical}%
  \fi \next}

% operatori

\def\multimathop #1 {\def\arg{#1}%
  \ifx\arg\pto \let\next\relax
  \else
  \def\next{\expandafter
    \def\csname #1\endcsname{\mathop{\rm #1}\nolimits}%
    \multimathop}%
  \fi \next}

\multibold
qwertyuiopasdfghjklzxcvbnmQWERTYUIOPASDFGHJKLZXCVBNM.

\multical
QWERTYUIOPASDFGHJKLZXCVBNM.

\multimathop
diag dist div dom mean meas sign supp .

% accorpamenti di formule citate:
% uso  \accorpa {prima}{seconda}
%      \Accorpa\cs prima seconda (con il comodo blank anche dopo)
% NB: \Accorpa definisce \cs come l'accorpamento delle due citazioni
% e scrive sul file.log

\def\accorpa #1#2{\eqref{#1}--\eqref{#2}}
\def\Accorpa #1#2 #3 {\gdef #1{\eqref{#2}--\eqref{#3}}%
  \wlog{}\wlog{\string #1 -> #2 - #3}\wlog{}}

% macro comode

\def\separa{\noalign{\allowbreak}}

\def\somma #1#2#3{\sum_{#1=#2}^{#3}}

\def\graffe #1{\mathopen\{#1\mathclose\}}

\def\<#1>{\mathopen\langle #1\mathclose\rangle}
\def\norma #1{\mathopen \| #1\mathclose \|}

\def\[#1]{\mathopen\langle\!\langle #1\mathclose\rangle\!\rangle}

\def\ioT {\int_0^T}

\def\iO{\int_\Omega}

\def\dt{\partial_t}
\def\dn{\partial_\nu}

\def\cpto{\,\cdot\,}

\def\checkmmode #1{\relax\ifmmode\hbox{#1}\else{#1}\fi}
\def\aeO{\checkmmode{a.e.\ in~$\Omega$}}
\def\aeQ{\checkmmode{a.e.\ in~$Q$}}

\def\aet{\checkmmode{a.e.\ in~$(0,T)$}}

\def\aat{\checkmmode{for a.a.~$t\in(0,T)$}}

% insiemi numerici

\def\erre{{\mathbb{R}}}

\def\enne{{\mathbb{N}}}

% spazi di funzioni a valori vettoriali su [0,T], [0,t], [0,s], [0,+\infty), [\delta,T]

% Come ricordare: in generale i simboli L H W  C da soli per gli spazi su (0,T)
% gli stessi raddoppiati per (0,+\infty)
% aggiunta di t o s al simbolo per (0,t) e (0,s)
% aggiunta di d al simbolo semplice o doppio per intervalli (\delta,T) e (\delta,+\infty)
% il simbolo C e i suoi derivati mettono le quadre anziche' le tonde

% Esempi   \L2V   \L\infty\Vp   \W{1,1}H   \C0H   \LL2V   \CC0\Vp   \Ld2V  \CCdH

\def\genspazio #1#2#3#4#5{#1^{#2}(#5,#4;#3)}
\def\spazio #1#2#3{\genspazio {#1}{#2}{#3}T0}

\def\L {\spazio L}
\def\H {\spazio H}
\def\W {\spazio W}

\def\C #1#2{C^{#1}([0,T];#2)}

% spazi di funzioni su \Omega, \Gamma, Q e \Sigma

\def\Lx #1{L^{#1}(\Omega)}
\def\Hx #1{H^{#1}(\Omega)}

\def\LQ #1{L^{#1}(Q)}

\def\Ldue{\Lx 2}
\def\Linfty{\Lx\infty}

\def\Huno{\Hx 1}
\def\Hdue{\Hx 2}
\def\Hunoz{{H^1_0(\Omega)}}

% spazi di funzioni su Q e S

\def\LQ #1{L^{#1}(Q)}

% lettere greche

\let\theta\vartheta

\let\phi\varphi

\let\hat\widehat

\let\TeXchi\chi                         % new \chi, exactly on the baseline
\newbox\chibox
\setbox0 \hbox{\mathsurround0pt $\TeXchi$}
\setbox\chibox \hbox{\raise\dp0 \box 0 }
\def\chi{\copy\chibox}

% quadratino di fine dimostrazione

% abbreviazioni specifiche del lavoro

\def\Az #1{A_0^{#1}}
\def\Vz #1{V_0^{#1}}
\def\VA #1{V_A^{#1}}
\def\VB #1{V_B^{#1}}
\def\VBz #1{V_{B,0}^{#1}}
\def\Hz{H_0}

\def\Pi{\hat\pi}

\def\yh{\hat \xi_{h_N}}

\def\overyh{\overline \xi_{h_N}}
\def\overmuh{\overline\eta_{h_N}}
\def\underyh{\underline \xi_{h_N}}
\def\undermuh{\underline\eta_{h_N}}
\def\overuh{\overline k_{h_N}}

\def\yn{\xi^n_N}
\def\ynp{\xi^{n+1}_N}

\def\mun{\eta^n_N}
\def\munp{\eta^{n+1}_N}

\def\dhyn{\frac{\ynp-\yn}{h_N}}

\def\yk{\xi^k}
\def\muk{\eta^k}

\def\yz{y_0}

\def\uad{{\cal U}_{\rm ad}}
\def\muk{\mu^k}
\def\yk{y^k}
\def\bu{\bar u}
\def\bm{\bar \mu}
\def\by{\bar y}
\def\zk{z^k}
\def\rk{\rho^k}
\def\etak{\eta^k}
\def\xik{\xi^k}
\def\pO{p_\Omega}
\def\1{{\bf 1}}
\def\calVp{\calV^*}

%%%%%%%%%%%%%%%%%%%%%%%%%%%%%%
\Begin{document}
%%%%%%%%%%%%%%%%%%%%%%%%%%%%%%%%%

%%%%%%%%%%%%%%%%%%%%%%%%%%%%%%%%%
%% front page
%%%%%%%%%%%%%%%%%%%%%%%%%%%%%%%%%

%
\title{Optimal distributed control of a generalized
 \\ fractional Cahn--Hilliard system}
\author{}
\date{}
\maketitle
\Bcenter
\vskip-1cm
{\large\sc Pierluigi Colli$^{(1)}$}\\
{\normalsize e-mail: {\tt pierluigi.colli@unipv.it}}\\[.25cm]
{\large\sc Gianni Gilardi$^{(1)}$}\\
{\normalsize e-mail: {\tt gianni.gilardi@unipv.it}}\\[.25cm]
{\large\sc J\"urgen Sprekels$^{(2)}$}\\
{\normalsize e-mail: {\tt sprekels@wias-berlin.de}}\\[.45cm]
$^{(1)}$
{\small Dipartimento di Matematica ``F. Casorati'', Universit\`a di Pavia}\\
{\small and Research Associate at the IMATI -- C.N.R. Pavia}\\
{\small via Ferrata 5, 27100 Pavia, Italy}\\[.2cm]
$^{(2)}$
{\small Department of Mathematics}\\
{\small Humboldt-Universit\"at zu Berlin}\\
{\small Unter den Linden 6, 10099 Berlin, Germany}\\[2mm]
{\small and}\\[2mm]
{\small Weierstrass Institute for Applied Analysis and Stochastics}\\
{\small Mohrenstrasse 39, 10117 Berlin, Germany}
\Ecenter
\Begin{abstract}\noindent
In the recent paper \pier{``Well-posedness and regularity for a generalized
fractional Cahn--Hilliard system'' by the same authors}, 
general well-posedness results have been 
established for a \pier{a class of evolutionary systems} of two equations
having the structure of a viscous Cahn--Hilliard system, in which nonlinearities of
double-well type occur. The operators  appearing in the
system equations are fractional versions in the spectral sense of general linear operators \,$A,B$\, having
compact resolvents, which are densely defined,  unbounded, selfadjoint, and monotone
in a Hilbert space of functions defined in a smooth domain. In this \giapie{work} we 
complement the results given in \giapie{the quoted paper} 
by studying a distributed control problem for
this evolutionary system. The main difficulty in the analysis is to establish a rigorous
Fr\'echet differentiability result for the associated control-to-state mapping. This
seems only to be possible if the state stays bounded, which, in turn, makes it necessary 
to postulate an additional global boundedness assumption.
One typical situation, in which this assumption is satisfied,
arises when $B$ is the negative Laplacian with
zero Dirichlet boundary conditions and the nonlinearity is smooth with
polynomial growth of at most order four. \pier{Also a case with logarithmic 
nonlinearity can be handled.} Under \pier{the global boundedness assumption}, 
we establish existence
and first-order necessary optimality conditions for the optimal control problem
in terms of a variational inequality and the associated adjoint state system.

\End{abstract}
\salta
\pagestyle{myheadings}
\newcommand\testopari{\sc Colli \ --- \ Gilardi \ --- \ Sprekels}
\newcommand\testodispari{\sc Optimal control of a fractional Cahn--Hilliard system}
\markboth{\testopari}{\testodispari}
\finqui
%
%%%%%%%%%%%%%%%%%%%%%%%%%%%%%%%%%
%% very beginning
%%%%%%%%%%%%%%%%%%%%%%%%%%%%%%%%%

\section{Introduction}
\label{Intro}
\setcounter{equation}{0}

Let $\Omega\subset \erre^3$ denote an open, bounded, and connected set with smooth
boundary $\,\Gamma\,$  and outward normal derivative $\,\dn$, let $T>0$  be a 
final time,  and let $H:=L^2(\Omega)$ denote the Hilbert space of square-integrable 
real-valued functions defined on $\Omega$, endowed with the standard inner product
$(\cdot,\cdot)$ and norm $\,\|\cdot\|$, respectively. We \pier{set} 
\,$Q_t:=\Omega\times (0,t)\,$ for $0<t<T$ and \,$Q:=\Omega\times (0,T)$. 
We investigate in this paper the following abstract distributed optimal control problem: 

\noindent{\bf (CP)}  Minimize the tracking-type cost functional
\revis{%
\Beq
\label{cost}
  {\cal J}((\mu,y),u)
  := \frac{\alpha_1}2\,\|y(T)-y_\Omega\|^2
  \,+\,\frac {\alpha_2}2\int_0^T\!\|y(t)-y_Q(t)\|^2\,dt
  \,+\,\frac{\revis{\alpha_ 3}}2\int_0^T\!\|u(t)\|^2\,dt
\Eeq
}%
over the admissible set
\begin{equation}
\label{uad}
\uad:=\left\{u\in H^1(0,T;\revis{H}):\,|u|\le\rho_1 \mbox{ \,a.\,e. in }\,Q, \,\,\,\,\|u\|_
{H^1(0,T;\revis{H})}\,\le\,\rho_2\right\},
\end{equation}
subject to the evolutionary state system
\begin{align}
  & \dt y + A^{2r} \mu = 0,
  \label{Iprima}
  \\
  & \tau \dt y + B^{2\sigma} y + f'(y) = \mu + u,
  \label{Iseconda}
  \\
  & y(0) = \yz.
  \label{Icauchy}
\end{align}
Here, \pier{$\rho_1$ and $\rho_2$ are fixed positive constants;} $\alpha_i$, \revis{$i=1,2,3$}, are 
nonnegative \pier{coefficients} but not all zero, and the given target functions satisfy $y_\Omega\in H$ and \revis{$y_Q\in L^2(0,T;H)$}. 
\revis{We remark at this point that the $\H1{\revis{H}}$ regularity of the admissible controls cannot be avoided due to analytical reasons.}
The linear operators $A^{2r}$ and $B^{2\sigma}$, 
with $r>0$ and $\sigma>0$, denote fractional powers (in the spectral sense) of operators $A$ and $B$. We will 
give a proper definition of such operators in the next section. Throughout this paper, we generally assume:
 
\vspace{1mm}\noindent
{\bf (A1)} \,\,\,$A:D(A)\subset H\to H$\, and \,$B:D(B)\subset H\to H\,$ are unbounded, monotone,\linebreak
\hspace*{12.3mm} and selfadjoint linear operators with compact resolvents.

\vspace{1mm}\noindent
This assumption implies that there are sequences 
$\{\lambda_j\}$ and $\{\lambda'_j\}$ of eigenvalues
and orthonormal sequences $\{e_j\}$ and $\{e'_j\}$ of corresponding eigenvectors,
that~is,
\Beq
  A e_j = \lambda_j e_j, \quad
  B e'_j = \lambda'_j e'_j,
  \aand
  (e_i,e_j) = (e'_i,e'_j) = \delta_{ij},
  \quad \hbox{for $i,j=1,2,\dots$,}
  \label{eigen}
\Eeq
such that
\begin{align}
  & 0 \leq \lambda_1 \leq \lambda_2 \leq \dots,
  \aand
  0 \leq \lambda'_1 \leq \lambda'_2 \leq \dots,
  \quad \hbox{with} \quad
  \lim_{j\to\infty} \lambda_j
  = \lim_{j\to\infty} \lambda'_j
  = + \infty,
  \label{eigenvalues}
  \\[1mm]
  & \hbox{$\{e_j\}$ and $\{e'_j\}$ are complete systems in $H$}.
  \label{complete}
\end{align}
 
Note that the state system \eqref{Iprima}--\eqref{Icauchy} 
can be seen as a generalization
of the famous Cahn--Hilliard system which models a  
phase separation process taking place in the container $\Omega$. In this case,
one typically has  $A^{2r}=B^{2\sigma}=-\Delta$ with zero Neumann or Dirichlet boundary conditions, and the unknown 
functions $\,y\,$ and $\,\mu\,$ stand for the \emph{order parameter} (usually a scaled
density of one of the involved phases) and the \emph{chemical potential}
associated with the phase transition, respectively. Moreover, $\,f\,$ denotes a 
double-well potential.
Typical and physically significant examples for $\,f\,$  
are the so-called {\em classical regular potential}, the {\em logarithmic double-well potential\/},
and the {\em double obstacle potential\/}, which are given, in this order,~by
\begin{align}
  & f_{reg}(r) := \frac 14 \, (r^2-1)^2 \,,
  \quad r \in \erre, 
  \label{regpot}
  \\
  & f_{log}(r) := \bigl( (1+r)\ln (1+r)+(1-r)\ln (1-r) \bigr) - c_1 r^2 \,,
  \quad r \in (-1,1),
  \label{logpot}
  \\[1mm]
  & f_{2obs}(r) := - c_2 r^2 
  \quad \hbox{if $|r|\leq1$}
  \aand
  f_{2obs}(r) := +\infty
  \quad \hbox{if $|r|>1$}.
  \label{obspot}
\end{align}
Here, the constants $c_1>1$ and $c_2>0$ in \eqref{logpot} and \eqref{obspot} are such 
that $f_{log}$ and $f_{2obs}$ are nonconvex. \juerg{Notice that in the case of the nondifferentiable potential
\eqref{obspot} the state equation \eqref{Iseconda} has to be understood as a variational inequality.}
We also note that $\tau$ is a nonnegative parameter, 
where for the classical
Cahn--Hilliard system one has $\tau=0$ (the \emph{nonviscous} case), while
$\tau>0$ corresponds to the \emph{viscous} case.

In the recent paper \cite{CGS18}, general well-posedness and regularity results
for the state system \eqref{Iprima}--\eqref{Icauchy} have been established for both the viscous and
nonviscous cases and for nonlinearities that include all of the three cases \eqref{regpot}--\eqref{obspot}.
It turned out that the first eigenvalue $\lambda_1$ of $\,A\,$ plays an important role
in the analysis. 
Indeed, the main assumption for the operators $A,B$ besides {\bf (A1)} was the following
\revis{(see Remark~\ref{CommentA2} below for examples)}:

\vspace{1mm}\noindent
{\bf (A2)} \,\,\,Either\\
\hspace*{20mm}(i) \,\,\,$\lambda_1>0$\\
\hspace*{12.5mm} or\\
\hspace*{20mm}(ii) \,\,$0=\lambda_1<\lambda_2$, and $\,e_1\,$ is a constant and belongs to the domain of $\,B^\sigma$. 

\vspace{1mm}
For our analysis of the optimal control problem {\bf (CP)}, the general assumptions {\bf (A1)} and {\bf (A2)}
are not sufficient. Indeed, in order to
be able to prove that the control-to-state operator ${\cal S}:u\mapsto (\mu,y)$ is Fr\'echet differentiable
between suitable Banach spaces, \pier{it} seems to be indispensable to assume that $\,f\,$ is smooth in its domain 
(which means that the potential \eqref{obspot} is not admitted) and 
to have at disposal an $L^\infty(Q)$ bound for both the state component $y$ and the functions $\,f^{(i)}(y)$,
for $i=1,2,3$. In the case of the logarithmic potential \eqref{logpot}, this means that we need to separate $\,y\,$ away
from the critical arguments $\pm 1$. We will discuss in Section 3 three situations in which \oldgianni{appropriate} boundedness  
conditions for $y$ and the derivatives $f'(y)$ can be guaranteed, where one of these cases applies to the logarithmic
potential.

Under these boundedness assumptions, we will be able to show
the Fr\'echet differentiability of the control-to-state operator ${\cal S}$ (cf.\ Section~\ref{FRECHET}) and to derive 
first-order necessary optimality conditions (cf.\ Section~\ref{OPTIMUM}).    

Let us add a few remarks on the existing literature. There exist numerous contributions on viscous/nonviscous,
local/nonlocal, convective/nonconvective Cahn--Hilliard systems for the classical (non-fractional) case 
$A=B=-\Delta$, $2r=2\sigma=1$, \pier{or some nonlocal counterparts,} where various types of boundary
conditions (e.g., Dirichlet, Neumann, dynamic) and different assumptions on the nonlinearity $f$ were considered. 
We refer the interested reader to\revis{\cite{ABG} and, for a selection of associated references,
to the recent paper\cite{CGSAnnali}}. 
Some papers also address the coupled Cahn--Hilliard/Navier--Stokes system 
(see, e.g,\pier{\cite{FGG, FGGS}} and the references given therein).

The literature on optimal control problems for non-fractional Cahn--Hilliard system is \revis{still scarce}.    
The case of Dirichlet and/or Neumann boundary conditions for various types of such systems were the
subject of, e.g., the works\pier{\cite{CGRS, CGSAIMS, CGSEECT, Duan, WN, Z, ZW}}, while the case of  
dynamic boundary conditions was studied in \cite{CFGS1, CFGS2, CFGS3, CGSAdvan, CGSAMO, CGSSIAM18, CGSConvex,
CS, Fukao}. The optimal control of convective Cahn--Hilliard systems was addressed in \cite{RS, ZL1, ZL2},
while the papers \cite{FGS, FRS, HHCK, HKW, HW1, HW2, HW3, Medjo} were concerned with 
coupled Cahn--Hilliard/Navier--Stokes systems. 

There are only a few contributions to the theory of Cahn--Hilliard systems involving fractional operators.
In the connection of well-posedness and regularity results, we refer to \cite{AM,AkSeSchi} for the case of 
the fractional negative Laplacian with zero Dirichlet boundary conditions; general operators other than the 
negative Laplacian
have apparently only studied in\revis{\cite{GalDCDS, GalEJAM, GalAIHP, CGS18}}. 
As of now, aspects of optimal control have been scarcely
dealt with even for simpler linear evolutionary systems involving fractional operators; for such systems,
some identification problems were addressed in the recent contributions \oldgianni{\cite{GV, PV, SV}}, while
for optimal control problems for such cases we refer to \cite{HAS} (for the stationary \pier{-- elliptic --} case, see
also \cite{HO1, HO2}).
However, to the authors' best
knowledge, the present paper appears to be the first contribution that addresses optimal control problems
for Cahn--Hilliard systems with general fractional order operators.

The paper is organized as follows: the subsequent Section~\ref{FPAM} brings some auxiliary functional analytic material, 
while in Section~\ref{STATE} some preparatory results concerning the state system \eqref{Iprima}--\eqref{Icauchy} are discussed. 
In Section~\ref{FRECHET}, the Fr\'echet differentiability of the control-to-state operator is shown, 
and in the final Section~\ref{OPTIMUM}, we then prove an existence result for the optimal control problem 
and establish the first-order necessary conditions of optimality.

\oldgianni{Throughout the paper, \pier{for a general Banach space $\,X\,$ we denote}
by $\,\|\cdot\|_X\,$ and  $\,X^*\,$ its norm and dual space, respectively. 
However, particular symbols are adopted for the spaces we introduce in the next section.}
% Finally, we use the notation $Q_t:=\Omega\times (0,t)$ 
% for $t\in (0,T]$ (where we have $Q=Q_T$).

%%%%%%%%%%%%%%%%%%%%%%%%%%%%%%%%%%%%%%%%%%%%%%%%%%%%%%%%%%%%%%%%%%%%%%%%

\section{Fractional powers and auxiliary results}
\label{FPAM}
\setcounter{equation}{0}

In this section, we collect some auxiliary material concerning functional analytic notions.
To this end, we generally assume that the conditions {\bf (A1)} and {\bf (A2)} are satisfied.
At this point, some remarks on the assumption {\bf (A2)} are in order.

\Brem
\label{CommentA2}
First, the meaning of {\bf (A2)},(i) is clear, and this condition is satisfied
for the more usual elliptic operators with zero Dirichlet boundary conditions
(however, also zero mixed boundary conditions could be considered, 
with proper definitions of the domains of the operators).
For instance, \,$A\,$ can be the Laplace operator $-\Delta$ with domain $D(-\Delta)=\Hdue\cap\Hunoz$ or
the bi-harmonic operator $\Delta^2$ with the domain\, 
$D(\Delta^2)=\Hx4\cap H^2_0(\Omega)$.
The second case {\bf (A2)},(ii),
in which the strict inequality means that the first eigenvalue $\lambda_1=0$ is simple,
arises in both of the following important situations:
$A\,$ is the Laplace operator $-\Delta$ with zero Neumann boundary conditions,
which corresponds to the choice $D(-\Delta)=\{v\in\Hdue:\ \dn v=0 \mbox{ \,on\, }\Gamma\}$, or
$\,A\,$ is the bi-harmonic operator $\Delta^2$ with the boundary conditions
encoded in the definition of the domain  
$D(\Delta^2)=\{v\in\Hx4:\ \dn v=\dn\Delta v=0\mbox{ \,on\, }\Gamma\}$.
Indeed, $\Omega$~is assumed to be bounded, smooth and connected.
\Erem

%\vspace{4mm} 
Using the facts summarized in \eqref{eigen}--\eqref{complete},
we can define the powers of $A$ and $B$ 
for an arbitrary positive real exponent.
For the first operator, we have
\begin{align}
  & \VA r := D(A^r)
  = \Bigl\{ v\in H:\ \somma j1\infty |\lambda_j^r (v,e_j)|^2 < +\infty \Bigr\}
  \pier{,}
  \label{defdomAr}
  \\
  & A^r v = \somma j1\infty \lambda_j^r (v,e_j) e_j
  \quad \hbox{for $v\in\VA r$},
  \label{defAr}
\end{align}
the series being convergent in the strong topology of~$H$,
due to the properties \eqref{defdomAr} of the coefficients.
In principle, we can endow $\VA r$ with the (graph) norm and inner product
\Beq
  \norma v_{gr,A,r}^2 := (v,v)_{gr,A,r}
  \aand
  (v,w)_{gr,A,r} := (v,w) + (A^r v , A^r w)
  \quad \hbox{for $v,w\in\VA r$}.
  \label{defnormagrAr}
\Eeq
This makes $\VA r$ a Hilbert space.
However, we can choose any equivalent Hilbert norm.
Indeed, in view of assumption {\bf (A2)}, it is more convenient to work with the Hilbert norm 
\Beq
  \norma v_{A,r}^2 := \left\{ 
  \begin{aligned}
  & \norma{A^r v}^2
  = \somma j1\infty |\lambda_j^r (v,e_j)|^2
  \qquad \hbox{if $\lambda_1>0$,}
  \\
  & |(v,e_1)|^2 + \norma{A^r v}^2
  = |(v,e_1)|^2 + \somma j2\infty |\lambda_j^r (v,e_j)|^2
  \qquad \hbox{if $\lambda_1=0$}.
  \end{aligned}
  \right.
  \label{defnormaAr}
\Eeq
In \cite[Prop.~3.1]{CGS18} it has been shown that this norm is equivalent 
to the graph norm defined in~\eqref{defnormagrAr},
and we always will work with the norm \eqref{defnormaAr} \pier{instead of} \eqref{defnormagrAr}.
We also use the corresponding inner product in $\VA r$ 
given~by
\Bsist
  && (v,w)_{A,r}
  = (A^r v,A^r w)
  \quad \hbox{or} \quad
  (v,w)_{A,r}
  = (v,e_1)(w,e_1) + (A^r v,A^r w),
  \qquad
  \non
  \\
  \label{inpro}
  &&  \hbox{depending on whether $\lambda_1>0$ or $\lambda_1=0$,\quad 
	for $v,w\in\VA r$.}
\Esist

\Brem
\label{Meanvalue}
Observe that in the case $\lambda_1=0$ 
the constant value of $e_1$
equals one of the numbers $\pm|\Omega|^{-1/2}$,
where $|\Omega|$ is the volume of~$\Omega$.
It follows for every $v\in H$ that the first term $(v,e_1)e_1$ of the Fourier series of $v$
is the constant function whose value~is the mean value of~$v$, which is defined by
\Beq
  \mean (v) := \frac 1 {|\Omega|} \iO v\,.
  \label{defmean}
\Eeq
Moreover, the first terms of the sums appearing in \eqref{defnormaAr} and \eqref{inpro}
are given~by
\Bsist
  && |(v,e_1)|^2 = |\Omega| \, (\mean v)^2
  \quad \hbox{for every $v\in H$},
  \non
  \\[1mm]
  && (v,e_1) (w,e_1) = |\Omega| \, (\mean v) (\mean w)
  \quad \hbox{for every $v,w\in H$}.
  \non
\Esist
\Erem

%\vspace{3mm}
In the same way as for $A$, starting from \accorpa{eigen}{complete} for~$B$,
we can define the power $B^\sigma$ of $B$ for every $\sigma>0$, where for 
$V_B^\sigma$ we choose the graph norm.
We therefore set
\Bsist
  && \VB\sigma := D(B^\sigma),
  \quad \hbox{with the norm $\norma\cpto_{B,\sigma}$ associated to the inner product}
  \label{defBs}
  \non
  \\
  && (v,w)_{B,\sigma} := (v,w) + (B^\sigma v,B^\sigma w)
  \quad \hbox{for $v,w\in \VB\sigma$}.
  \label{defprodBs}
\Esist

\Brem
Let us briefly comment on the condition {\bf (A2)},(ii).  
We notice that the condition that $e_1$ be a  constant belonging to $V_B^\sigma$ holds true  
for many operators having a domain \revis{involving} Neumann boundary conditions.
This is the case, for instance,
if $B$ is the Laplace operator with domain $D(-\Delta)=\{v\in\Hdue:\ \dn v=0 \mbox{\, on \,}\Gamma\}$.
On the contrary, if $B=-\Delta$ with domain $D(-\Delta):=\Hdue\cap\Hunoz$,
then $D(B)$ does not contain any nonzero constant functions.
However, $\VB\sigma$ does contain every constant function
provided that $\sigma\in(0,1/4)$, since $\VB\sigma$ is in this case a subspace
of the usual Sobolev-Slobodeckij space~$\Hx{2\sigma}$.
\Erem

%\vspace{2mm}
To resume our preparations, we observe that if $r_i$ and $\sigma_i$ are arbitrary positive exponents,
then it is easily seen that we have the ``Green type'' formulas
\Bsist
  && (A^{r_1+r_2} v,w)
  = (A^{r_1} v, A^{r_2} w)
  \quad \hbox{for every $v\in\VA{r_1+r_2}$ and $w\in\VA{r_2}$},
  \label{propA}
  \\[1mm]
  && (B^{\sigma_1+\sigma_2} v,w)
  = (B^{\sigma_1} v, B^{\sigma_2} w)
  \quad \hbox{for every  $v\in\VB{\sigma_1+\sigma_2}$ and $w\in\VB{\sigma_2}$}.
  \label{propB}
\Esist

The next step is the introduction of some spaces with negative exponents.
We set
\Beq
  \VA{-r} := (\VA r)^* 
  \quad \hbox{for $r>0$},
  \label{defVAneg}
\Eeq
and endow $\VA{-r}$ with the dual norm $\,\|\cdot\|_{A,-r}\,$ of $\,\|\cdot\|_{A,r}$.
\oldgianni{%
We use the symbol $\<\cpto,\cpto>_{A,r}$ for the duality pairing
between $\VA{-r}$ and~$\VA r$ and identify $H$ with a subspace of $\VA{-r}$
in the usual sense, i.e., in order that 
$\<z,v>_{A,r}=(z,v)$ for every $z\in H$ and $v\in\VB\sigma$.
Similarly, we~set
\Beq
  \VB{-\sigma} := (\VB\sigma)^* 
  \quad \hbox{for $\sigma>0$}.
  \label{defVBneg}
\Eeq
As $\VB\sigma$ is dense in~$H$, we have the analogous embedding}
\Beq
  \oldgianni{H \subset \VB{-\sigma} .
  \label{embVBneg}}
\Eeq

Observe that the following embedding results are valid:
\begin{align}
  & \hbox{the embeddings $\VA{r_2} \subset \VA{r_1} \subset H$ are dense and compact for $0<r_1<r_2$}\pier{;}
  \label{compembA}
  \\
  & \hbox{the embeddings $H \subset \VA{-r_1} \subset \VA{-r_2}$ are dense and compact for $0<r_1<r_2$}\pier{;}
  \qquad
  \label{compembAneg}
  \\
  & \hbox{the embeddings $\VB{\sigma_2} \subset \VB{\sigma_1} \subset H$ are dense and compact for $0<\sigma_1<\sigma_2$}.
  \label{compembB}
\end{align}
%
%Moreover, from the continuous embedding $H\subset\VA{-r}$ 
%and the compact embedding $\VB\sigma\subset H$ given by~\accorpa{compembAneg}{compembB},
%it follows that, for every $\delta>0$, there exists a constant $c_\delta$ such that
%\Beq
%  \norma v^2
%  \leq \delta \, \norma{B^\sigma v}^2 
%  + c_\delta \norma v_{\VA{-r}}^2
%  \quad \hbox{for every $v\in\VB\sigma$}.
%  \label{compineq}
%\Eeq
We also note the validity of the Poincar\'e type inequality \oldgianni{(see \cite[formula~(3.5)]{CGS18})}
\Beq
  \norma v \,\leq\, \hat c \, \norma{A^r v}
  \quad \hbox{for every $v\in\VA r$ with $\mean (v)=0$}.
  \label{poincare} 
\Eeq
At this point, we introduce the Riesz isomorphism $\calR_r:\VA r\to\VA{-r}$
associated with the inner product~\eqref{inpro}, which is given by
\Beq
  \< \calR_r v , w >_{A,r}
  = (v,w)_{A,r}
  \quad \hbox{for every $v,w\in\VA r$}.
  \label{riesz}
\Eeq
Moreover, we~set
\oldgianni{%
\Bsist
  && \Vz r := \VA r
  \aand
  \Vz{-r} := \VA{-r}
  \quad \hbox{if $\lambda_1>0$},
  \label{defVrpos}
  \\[1mm]
  && \Vz r := \{v\in \VA r :\ \mean (v)=0\} \,\,\mbox{ and }
  \non
  \\
  &&   \Vz{-r} := \{v \in \VA{-r} :\ \<v,1>_{A,r}=0 \}
  \quad \hbox{if $\lambda_1=0$} \,.
  \label{defVrzero}
\Esist
}%
According to \cite[Prop.~3.2]{CGS18},  $\calR_r$ maps $\Vz r$ onto $\Vz{-r}$
and extends to $\Vz r$ the restriction of $A^{2r}$ to $\Vz{2r}$.
In view of this result, it is reasonable to use a proper notation
for the restrictions of $\calR_r$ and $\calR_r^{-1}$ to the subspaces
$\Vz r$ and $\Vz{-r}$, respectively.
We~set 
\Beq
  \Az{2r} := (\calR_r)_{|\Vz r}
  \aand
  \Az{-2r} := (\calR_r^{-1})_{|\Vz{-r}}\,,
  \label{defAz}
\Eeq
where the index $0$ has no meaning if $\lambda_1>0$ (since then $\Vz{\pm r}=\VA{\pm r}$),
while it reflects the zero mean value condition in the case $\lambda_1=0$.
We thus have
\Bsist
  && \Az{2r} \in \calL(\Vz r,\Vz{-r}) , \quad
  \Az{-2r} \in \calL(\Vz{-r},\Vz r)
  \aand
  \Az{-2r} = (\Az{2r})^{-1}\,,
  \label{contlinAz}
  \\[1mm]
  && \< \Az{2r} v,w >_{A,r} = (v,w)_{A,r} = (A^r v,A^r w)
  \quad \hbox{for every $v\in\Vz r$ and $w\in\VA r$}\,,
  \qquad
  \label{identityAz}
  \\[1mm]
  && \< f , \Az{-2r} f >_{A,r}
  = \norma{\Az{-2r} f}_{A,r}^2
  = \norma f_{A,-r}^2
  \quad \hbox{for every $f\in\Vz{-r}$}.
  \label{normaAz}
\Esist
Notice that \eqref{normaAz} implies that
\Beq
  \< f' , \Az{-2r} f >_{A,r}
  = \frac 12 \, \frac d {dt} \,\norma f_{A,-r}^2
  \quad \hbox{\aet, \ for every $f\in\H1{\Vz{-r}}$}.
  \label{dtnormaAz}
\Eeq
Moreover, by virtue of \cite[Prop.~3.3]{CGS18}, we have
\Beq
  \bigl( A^r \Az{-2r} f , A^r v )
  = \< f,v >_{A,r}
  \quad \hbox{for every $f\in\Vz{-r}$ and $v\in\VA r$}.
  \label{ok}
\Eeq
In addition (see \cite[Prop.~3.4]{CGS18}\oldgianni), the operator $A^{2r}\in\calL(\VA{2r},H)$ can be extended in a unique way
to a continuous linear operator, still termed~$A^{2r}$, from $\VA r$ into~$\Vz{-r}$, and we have
\Beq
  \norma{A^{2r}v}_{A,-r} \leq \norma{A^r v}
  \quad \hbox{for every $v\in\VA r$}.
  \label{stimaAdr}
\Eeq

As a final preparation, we now introduce some notations concerning interpolating functions.

\step 
Interpolants

Let \,$N$\, be a positive integer and $\,Z\,$ be one of the spaces \,$H$, $\VA r$,~$\VB\sigma$.
We set $\,h_N:=T/N\,$ and $\,I_N^n:=((n-1)h_N,nh_N)\,$ for $\,n=1,\dots,N$.
Given $\,z=(z_0,z_1,\dots ,z_N)\in Z^{N+1}$,
we define the piecewise constant and piecewise linear interpolants
\Beq
  \overline z_{h_N} \in \L\infty Z , \quad
  \underline z_{h_N} \in \L\infty Z 
  \aand
  \hat z_{h_N} \in \W{1,\infty}Z
  \non
\Eeq
by setting 
\Bsist
  && \hskip -2em
  \overline z_{h_N}(t) = z^n
  \aand
  \underline z_{h_N}(t) = z^{n-1}
  \quad \hbox{for a.a.\ $t\in I_N^n$, \ $n=1,\dots,N$},
  \label{pwconstant}
  \\
  && \hskip -2em
  \hat z_{h_N}(0) = z_0
  \aand
  \dt\hat z_{h_N}(t) = \frac {z^{n+1}-z^n}{h_N}
  \quad \hbox{for a.a.\ $t\in I_N^n$, \ $n=1,\dots,N$}.
  \qquad
  \label{pwlinear}
\Esist

For the reader's convenience,
we summarize some well-known relations between the finite set of values
and the interpolants. We have that
\begin{align}
  & \norma{\overline z_{h_N}}_{\L\infty Z}
  = \max_{n=1,\dots,N} \norma{z^n}_Z \,, \quad
  \norma{\underline z_{h_N}}_{\L\infty Z}
   = \max_{n=0,\dots,N-1} \norma{z^n}_Z\,,
  \label{ouLinftyZ}
  \\
  & \norma{\dt\hat z_{h_N}}_{\L\infty Z}
  = \max_{0\leq n\leq N-1} \norma{(z^{n+1}-z^n)/h_N}_Z\,,
  \label{dtzLinftyZ}
  \\
  \separa
  & \norma{\overline z_{h_N}}_{\L2Z}^2
  = h_N \somma n1N \norma{z^n}_Z^2 \,, \quad
  \norma{\underline z_{h_N}}_{\L2Z}^2
  = h_N \somma n0{N-1} \norma{z^n}_Z^2 \,,
  \label{ouLdueZ}
  \\
  \separa
  & \norma{\dt\hat z_{h_N}}_{\L2Z}^2
  = h_N \somma n0{N-1} \norma{(z^{n+1}-z^n)/h_N}_Z^2\,, 
  \label{dtzLdueZ}
  \\
  & \norma{\hat z_{h_N}}_{\L\infty Z}
  = \max_{n=1,\dots,N} \max\{\norma{z^{n-1}}_Z,\norma{z^n}_Z\}
  = \max\{\norma{z_0}_Z,\norma{\overline z_{h_N}}_{\L\infty Z}\}\,,
  \label{hzLinftyZ}
  \\
  & \norma{\hat z_{h_N}}_{\L2Z}^2
  \leq h_N \somma n1N \bigl( \norma{z^{n-1}}_Z^2 + \norma{z^n}_Z^2 \bigr)
  \leq h_N \norma{z_0}_Z^2
  + 2 \norma{\overline z_{h_N}}_{\L2Z}^2 \,.
  \label{hzLdueZ}
\end{align}
Moreover, it holds that
\begin{align}
  & \norma{\overline z_{h_N}-\hat z_{h_N}}_{\L\infty Z}
  = \max_{n=0,\dots,N-1} \norma{z^{n+1}-{z^n}}_Z
  = h_N \, \norma{\dt\hat z_{h_N}}_{\L\infty Z}\,,
  \label{diffLinfty}
  \\
  & \norma{\overline z_{h_N}-\hat z_{h_N}}_{\L2Z}^2
  = \frac {h_N} 3 \somma n0{N-1} \norma{z^{n+1}-z^n}_Z^2
  = \frac {h^2_N} 3 \, \norma{\dt\hat z_{h_N}}_{\L2Z}^2\,,
  \label{diffLdue}
\end{align}
and similar identities for the difference $\underline z_{h_N}-\hat z_{h_N}$.
As a consequence, we also have the inequalities
\Bsist
  && \norma{\overline z_{h_N}-\underline z_{h_N}}_{\L\infty Z}
  \leq 2h_N \, \norma{\dt\hat z_{h_N}}_{\L\infty Z}\,,
  \qquad
  \label{diffbisLinfty}
  \\
  && \norma{\overline z_{h_N}-\underline z_{h_N}}_{\L2Z}^2
  \leq \frac {\oldgianni 4 h_N^2} 3 \, \norma{\dt\hat z_{h_N}}_{\L2Z}^2 \,.
  \label{diffbisLdue}
\Esist
Finally, \pier{we observe that}
\Bsist
  && h_N \somma n0{N-1} \norma{(z^{n+1}-z^n)/h_N}_Z^2 
  \leq \norma{\dt z}_{\L2Z}^2
  \non
  \\
  && \quad \hbox{if $z\in\H1Z$\aand $z^n=z(nh_N)$ \quad for \, $n=0,\dots,N$}.
  \label{interpH1Z}
\Esist

\vspace{1mm}
Throughout the paper, we make use of
the elementary identity and inequalities
\Bsist
  \hskip-1cm&& a (a-b)
  = \frac 12 \, a^2
  + \frac 12 \, (a-b)^2
  - \frac 12 \, b^2
  \geq \frac 12 \, a^2
  - \frac 12 \, b^2
  \quad \hbox{for every $a,b\in\erre$},
  \label{elementare}
  \\
  \hskip-1cm&& ab\leq \delta a^2 + \frac 1 {4\delta}\,b^2
  \quad \hbox{for every $a,b\in\erre$ and $\delta>0$},
  \label{young}
\Esist
\Accorpa\Elementari elementare young
and quote \eqref{young} as the Young inequality.
We also take advantage of the summation by parts formula
\Beq
  \somma n0{k-1} a_{n+1} (b_{n+1} - b_n)
  = a_k b_k - a_1 b_0
  - \somma n1{k-1} (a_{n+1} - a_n) b_n\,,
  \label{byparts}
\Eeq
which is valid for arbitrary real numbers $a_1,\dots,a_k$ and $b_0,\dots,b_k$.
We also account for the discrete Gronwall lemma in the following form
(see, e.g., \cite[Prop.~2.2.1]{Jerome}):
for nonnegative real numbers $M$ and $a_n,b_n$, $n=0,\dots,N$, 
\Bsist
  a_k \leq M + \somma n0{k-1} b_n a_n 
  \quad \hbox{for $k=0,\dots,N$}
  \qquad \hbox{implies}
  \non
  \\
  a_k \leq M \exp \Bigl( \somma n0{k-1} b_n \Bigr)
  \quad \hbox{for $k=0,\dots,N$}.
  \label{gronwall}
\Esist
In \accorpa{byparts}{gronwall} it is understood that
a sum vanishes if the corresponding set of indices is empty.

\section{General assumptions and the state system}
\label{STATE}
\setcounter{equation}{0}
In this section, we state our general assumptions and discuss the
properties of the state system \eqref{Iprima}--\eqref{Icauchy}.
Besides {\bf (A1)} and {\bf (A2)}, we generally assume for the data of the state system:

\vspace{2mm}\noindent
{\bf (A3)}  \,\,\,$r$, $\sigma$, and $\tau$ are fixed positive real numbers.

\vspace{2mm}\noindent
\oldgianni{%
{\bf (A4)} \,\,\,$f=f_1+f_2$,\,\, where \,$f_1$, $f_2$\, and \,$f$\, satisfy:\\[1mm]
\hspace*{12mm}$\,f_1\in C^3(D(f_1))$, \,\,$D(f_1)$ being an open interval,\,\, and $\,f_1''\ge 0\,$ in \,$D(f_1)$;\\[0.5mm]
\hspace*{12mm}$\,f_2\in C^3(\erre)$, \,and\, $f_2'\,$ is Lipschitz continuous on $\erre$;\\[0.5mm]
\hspace*{12mm}$\,\,\liminf_{|s|\nearrow +\infty}\,\frac{f(s)}{s^2}>0$.
}%

\vspace{2mm}\noindent
{\bf (A5)} \,\,\,$y_0\in V^{2\sigma}_B\,$ and\, $f'(y_0)\in H$.

\vspace{3mm}\noindent
\oldgianni{Notice that {\bf (A4)} holds true for the classical regular potential \eqref{regpot},
for which we have $D(f_1)=\erre$.
In general, if $D(f_1)\not=\erre$, then 
it is understood that $f_1$ also stands for its l.s.c.\ extension in the sum $f=f_1+f_2$.
This is the case for the logarithmic potential \eqref{logpot}, 
for which we have $D(f_1)=(-1,1)$,
and its l.s.c.\ extension is given by setting 
$f_1(\pm 1):= 2\ln(2)$ and $f_1(r):=+\infty$ if $|r|>1$.
In cases like this, the growth condition at infinity for $f$ is trivially satisfied.
Finally, we remark that assumption {\bf (A4)} excludes the double obstacle potential~\eqref{obspot},
whose effective domain is the closed interval~$[-1,1]$.
}%

\vspace{2mm}
For the quantities entering the cost functional 
and the admissible set~$\uad$ \oldgianni{(see \eqref{cost} and~\eqref{uad})}, we generally assume:

\vspace{2mm}\noindent
{\bf (A6)} \,\,\,$y_\Omega\in\Ldue$, \,\revis{$y_Q\in L^2(Q)$}, \,the constants
$\,\alpha_i\ge 0$, \revis{$i=1,2,3$}, are not all equal \hspace*{12.4mm} to zero, \,$\rho_1>0$, and \,$\rho_2>0$.

\vspace{2mm}
Finally, we denote the control space by
\begin{equation}
{\cal X}\,:= \,H^1(0,T;\Ldue)\cap L^\infty(Q),
\label{defX}
\end{equation}
and make an assumption which is rather a denotation, since $\uad$ is a bounded
subset of~${\cal X}$:

\vspace{2mm}\noindent
{\bf (A7)}  \,\,\,The constant \,$R>0$ \,is such that \,$\uad\subset {\cal U}_R:=
\{u\in {\cal X}:\,\|u\|_{{\cal X}}\,<\,R\}.$

\vspace{3mm}
With the above assumptions, we are now ready to \revis{quote} a well-posedness result for the
state system \eqref{Iprima}--\eqref{Icauchy} which is a special case of the general
results \cite[Thm.~2.6 and Thm.~2.8]{CGS18}. \juerg{To this end, we recall the weak notion of solution 
to the system \eqref{Iprima}--\eqref{Icauchy} introduced in \cite{CGS18}. 
Namely, we look for a pair of functions $(\mu,y)$ satisfying the variational (in)equalities}
\begin{align}
\label{weak1}
&(\dt y(t),v)+(A^r\mu(t),\oldgianni{A^r}v)=0
\quad\mbox{for a.e. $\,t\in (0,T)\,$ and every $\,v\in V_A^r$},\\[2mm]
\label{weak2}
&\juerg{(\tau\dt y(t),y(t)-v)+(B^\sigma y(t),B^\sigma (y(t)-v))
+\giapie{\iO f_1(y(t))}\,+(f_2'(y(t)),y(t)-v)}\non\\[1mm]
&\juerg{\le\,(\mu(t)+u(t),y(t)-v)+\giapie{\iO f_1(v)}} \quad\mbox{for a.e. $\,t\in (0,T)\,$ and every $\,v\in V_B^\sigma$},\\[2mm]
\label{weak3}
&y(0)=y_0\,.
\end{align}

\vspace{2mm}
We have the following result. 

\Bthm
\label{Wellposedness}
Suppose that the general assumptions {\bf (A1)}--{\bf (A5)} and {\bf (A7)} are fulfilled. 
Then the weak state system
\eqref{weak1}--\eqref{weak3} has for every $\,u\in{\cal U}_R\,$ \revis{at least one} solution
$(\mu,y)$ such that
\begin{align}
\label{regmu}
&\mu\in L^\infty(0,T;V_A^{2r}),\\
\label{regy}
&y\in W^{1,\infty}(0,T;H)\cap H^1(0,T;V_B^\sigma),\\
\label{yL1}
&\juerg{f_1(y)\in L^1(Q)}.
\end{align}
Moreover, there are constants $\,K_1>0\,$ and $K_2>0$, which depend only on the data of the state system and $R$,
such that the following holds true:

\vspace{1mm}\noindent
{\rm (i)} \,\,Whenever $\,u\in{\cal U}_R\,$ is given, then the \revis{abovementioned} solution $(\mu,y)$ satisfies
\begin{align}
\label{bounds1}
\|\mu\|_{L^\infty(0,T;V_A^{2r})}\,+\,\|y\|_{W^{1,\infty}(0,T;H)\cap H^1(0,T;V_B^\sigma)}\,\le\,K_1.
\end{align}
\vspace{1mm}%
\noindent
{\rm (ii)} \,Whenever $\,u_i\in{\cal U}_R$, $i=1,2$, are given and $(\mu_i,y_i)$, $i=1,2$, \revis{are associated} solutions,
\phantom{{\rm (ii)} \,}then 
\begin{align}
\label{stabu1}
\|y_1-y_2\|_{L^\infty(0,T;H)\cap L^2(0,T;V_B^\sigma)}\,\le\,K_2\,\|u_1-u_2\|_{L^2(0,T;H)}\,.
\end{align}
\Ethm

\Brem \label{Pier-rem}
\pier{Note that the regularity \eqref{yL1} can be improved up to 
$$f_1(y)\in L^\infty(0,T;L^1(\Omega)).$$
Indeed,  \giapie{first, $f_1$~is bounded from below by an affine function, so that
$$ \iO f_1 (y(t) ) \geq - c \iO \bigl( 1+ |y(t)| \bigr) $$
for some constant $c>0$ and \revis{almost every} $t\in (0,T)$, and the last term is bounded since $y\in \L\infty H.$
On the other hand, thanks to \eqref{weak2}, $ \iO f_1 (y)$ is bounded} 
from above by an $L^\infty (0,T) $-function (cf.~\eqref{regmu}--\eqref{regy}).}
\Erem

The following global boundedness condition is crucial for the analysis of the control problem.

\revis{%
\bgroup
\setbox 0\hbox{{\bf (GB)} \,\,\,}%
\vspace{2mm}\noindent
\copy 0
\vtop{\advance\hsize -\wd0 \noindent
There is a compact interval $[a,b] \revis{\subset D(f_1)}$, 
which depends only on the data of the state system and~$R$, such that: 
if $(\mu,y)$ solves the state system for some $u\in{\cal U}_R$ in the sense of Theorem~\ref{Wellposedness},
then $y\in[a,b]$ \aeQ.}
\egroup
}%

\revis{%
\Brem
\label{OldGB}
If {\bf (GB)} holds, then there is some constant $K_3>0$ such that
\begin{align}
\label{bounds2}
\max_{i=0,1,2,3}\,\|f_1^{(i)}(y)\|_{L^\infty(Q)}\,\le\,K_3
\end{align}
for any solution $(\mu,y)$ in the sense of Theorem~\ref{Wellposedness} corresponding to some $u\in{\cal U}_R$.
\Erem
}%

\revis{%
\Brem
\label{CtoSmap}
If {\bf (GB)} holds, then we may argue as in \cite[Rem.~4.1]{CGS18} to conclude that
the whole pair $(\mu,y)$ is uniquely determined.
Therefore, the control-to-state operator
\Beq 
  {\cal S} : u \mapsto {\cal S}(u) := (\mu,y)
  \label{defS}
\Eeq
is well defined as a mapping from $\,{\cal U}_R\subset {\cal X}\,$ 
into the Banach space specified by the regularity conditions \accorpa{regmu}{regy}.
\Erem
}%

\revis{%
\Brem
\label{Addreg}
If {\bf (GB)} holds and
\vskip2mm
\noindent
{\bf (A8)} \,\,\,$\VB\sigma\cap\Linfty$ \ is dense in \ $\VB\sigma$
\vskip2mm
\noindent
then the variational inequality \eqref{weak2} implies the variational equality
\begin{align}
&(\tau\dt y(t),v)+(B^\sigma y(t),B^\sigma (v))+(f_1'(y(t)),v)+(f_2'(y(t)),v)
\,=\,(\mu(t)+u(t),v)\non\\
\label{weak2new}
&\quad\mbox{for a.e. $\,t\in (0,T)\,$ and every $\,v\in V_B^\sigma$},
\end{align}
as we show at once.
This implies that, a~fortiori, 
by virtue of \eqref{bounds1}, \eqref{bounds2} and a comparison in equation~\eqref{weak2new}, 
we have
$B^{2\sigma}y=\mu+u-\tau\,\dt y-f'(y)\in L^\infty(0,T;H)$, 
whence we can infer the additional regularity 
\Beq
  y \in \L\infty{V_B^{2\sigma}}.
  \label{addreg}
\Eeq
In particular, the solution $(\mu,y)$ is strong 
and \eqref{Iseconda} is valid almost everywhere in $Q$.
\hfil\break\indent
Let us deduce \eqref{weak2new} from~\eqref{weak2} by \revis{using}~{\bf(GB)}.
In the following, $t$~is fixed (\aet).
We write \eqref{weak2} with $w$ instead of~$v$.
Now, we pick any $v\in\VB\sigma\cap\Linfty$, term $M$ the $\Linfty$ norm of~$v$
and choose $\delta>0$ such that the compact interval $[a_\delta,b_\delta]:=[a-\delta M,b+\delta M]$
is contained in the open interval~$D(f_1)$.
At this point, we choose $w=y-\delta v$ in~\eqref{weak2} (so that $y-w=\delta v$)
and divide by~$\delta$.
We obtain
\Beq
  (\tau\dt y , v)
  + (B^\sigma y , B^\sigma v)
  - \iO \frac {f_1(y-\delta v) - f_1(y)} \delta
  + (f_2'(y) , v)
  \leq (\mu+u , v).
  \non
\Eeq
On the other hand, we have \aeO
\Beq
  \left| \frac {f_1(y-\delta v) - f_1(y)} \delta \right|
  \leq |v| \, \sup\{ f_1'(s):\ s\in[a_\delta,b_\delta \revis{]}\}.
  \non
\Eeq
By Lebesgue's theorem, we deduce that
\Beq
  \lim_{\delta\to0} \iO \frac {f_1(y-\delta v) - f_1(y)} \delta
  = - \iO f_1'(y) \, v
  = - (f_1'(y) , v).
  \non
\Eeq
Therefore, we conclude that
\Beq
  (\tau\dt y , v)
  + (B^\sigma y , B^\sigma v)
  + (f'(y) , v)
  \leq (\mu+u , v).
  \non
\Eeq
By applying this to $\pm v$ we \revis{infer} that
\Beq
  (\tau\dt y , v)
  + (B^\sigma y , B^\sigma v)
  + (f'(y) , v)
  = (\mu+u , v).
  \non
\Eeq
This means \eqref{weak2new}, but only for every $v\in\VB\sigma\cap\Linfty$, in principle.
However, by applying~{\bf(A8)}, we conclude that the same \revis{holds} for every $v\in\VB\sigma$.
\Erem
}%

\revis{%
\Brem
\label{RemA8}
Assumption {\bf(A8)} is satisfied in a variety of concrete situations.
An abstract sufficient condition is the following:
there exists $\sigma'$ such that $\VB{\sigma'}\subset\Linfty$.
Indeed, {\bf(A8)} trivially follows if $\sigma'\leq\sigma$ and also follows in the opposite case since
$\VB{\sigma'}$ is dense in~$\VB\sigma$ (cf.~\eqref{compembB}).
Another sufficient conditions is that the eigenfunctions of $B$ are bounded,
since the finite sums of the Fourier series of any $v\in\VB\sigma$
(which converge to $v$ in~$\VB\sigma$) 
are also bounded in this case.
\Erem
}%

\step
Examples

The condition {\bf (GB)} seems to be very restrictive and requires a case-to-case analysis.
 We now give \pier{some} sufficient conditions
under which it holds true. 
In all of the following three examples, we have $B=-\Delta$ with either zero Dirichlet or zero Neumann boundary condition.
Then\pier{, it turns out that} $V_B^1\subset H^2(\Omega)$, and thus, by regularity, $V_B^\sigma\subset H^{2\sigma}(\Omega)$ for all $\sigma\in\enne$.
Interpolation shows that then also $V_B^\sigma\subset H^{2\sigma}(\Omega)$ for all noninteger $\sigma>0$.
\revis{In particular, we have the embedding $\VB\sigma\subset\Linfty$ whenever $\sigma>3/4$, so that {\bf(A8)} trivially holds in this case
and the solution $(\mu,y)$ solves the variational equation~\eqref{weak2new} whenever {\bf(GB)} is satisfied.}
\juerg{We also notice that $V_{-\Delta}^{1/2}$ \pier{is equal to $H^1_0(\Omega)$ for Dirichlet boundary conditions or $\Huno$ in the case of Neumann boundary conditions.}}

\vspace{1mm}
{\bf 1.} \,\,We begin with the logarithmic potential \eqref{logpot}. \juerg{Recall that in this case we have
$f_1(r)=(1+r)\ln(1+r)+(1-r)\ln(1-r)$ for $r\in (-1,1)$, $f_1(\pm 1)=2\ln(2)$, and $f_1(r)=+\infty$ if
$r\not\in [-1,1]$. Hence it follows from 
the variational inequality \eqref{weak2} that the corresponding solution component $\,y\,$ must satisfy $\,y\in [-1,1]\,$ almost everywhere. In particular, $\|f_2'(y)\|_{L^\infty(Q)}$ is bounded.} Now assume
that $B=-\Delta$ with zero Neumann boundary condition, $2\sigma=1$, and 
\begin{equation}
 -1<\inf_{x\in\Omega} y_0(x), \quad \sup_{x\in\Omega} y_0(x)<+1.
\end{equation}
Moreover, assume that the embedding 
\begin{equation}
\label{ex1}
V_A^{2r}\subset L^\infty(\Omega)
\end{equation}
holds true. This is the case, for instance, if $A=-\Delta$ with zero Dirichlet or Neumann condition and $r>3/8$. Indeed,
we then have (see above) $V_A^{2r}\subset H^{4r}(\Omega)$ and $4r>3/2$,
\oldgianni{which implies that \,$\Hx{4r}\subset\Linfty$}. 
Now let \revis{$(\mu,y)$ be a solution corresponding to} some $u\in{\cal U}_R$. 
If \eqref{ex1} is satisfied, then we can infer from
\eqref{bounds1} that there is some global constant $M>0$ such that \juerg{$\,\|\mu+u-f_2'(y)\|_{L^\infty(Q)}\,\le\,M$.
By the form of the derivative $f_1'$} of the logarithmic potential, there are constants $r_*,r^*\in (-1,1)$ 
with $r_*\le y_0\le r^*$ a.e. in
$\Omega$ such that 
\begin{align*}
\juerg{f_1'(r)+M\le 0 \quad\forall\,r\in (-1,r_*)\quad\mbox{ and }\quad f_1'(r)-M\ge 0 \quad\forall\,r\in (r^*,1).}
\end{align*}
Now, \juerg{recall that $V_{-\Delta}^{1/2}= \Huno$. We thus may insert $v=y(t)-(y(t)-r^*)^+\in \Huno$ in the variational inequality  \eqref{weak2},
where  $(y(t)-r^*)^+$ is the positive part of $y(t)-r^*$. We then find for almost every $t\in (0,T)$ the inequality}
\begin{align}\label{FabFour}
&\juerg{\frac \tau 2\,\frac d{dt} \|(y(t)-r^*)^+\|^2\,+\,\iO |\giapie{\nabla(y(t)-r^*)^+}|^2}\non\\
&\juerg{\le\iO \bigl[ f_1(y(t)-(y(t)-r^*)^+)\,-\,f_1(y(t))\,+\,(\mu(t)+u(t)-f_2'(y(t)))(y(t)-r^*)^+\bigr]\,.}
\end{align}
\juerg{We claim that the integrand of the integral on the \rhs\ is nonpositive. To this end, we put
\begin{equation*} \Omega_+(t):=\{x\in\Omega:\,y(x,t)>r^*\}, \quad \Omega_-(t)=\{x\in\Omega:\,y(x,t)\le r^*\}.
\end{equation*}
Obviously, $(y(t)-r^*)^+=0$ on $\Omega_-(t)$, and thus the integrand is zero on $\Omega_-(t)$. On the other hand,
in $\Omega_+(t)$ we have $(y(t)-r^*)^+=y(t)- \pier{r^*}$, and thus the integrand equals}
\begin{equation*}
\juerg{f_1(r^*)-f_1(y(t))+(\mu(t)+u(t)-f_2'(y(t)))(y(t)-r^*).}
\end{equation*}
Now $r^*\in (-1,1)$, and thus $f_1$ is differentiable at $r^*$. Hence, invoking the convexity of $f_1$,
we have in $\Omega_+(t)$ that $\,f_1(r^*)-f_1(y(t))\,\le\,-f_1'(r^*)(y(t)-r^*)$. Now, by construction,
it holds that $\,\mu(t)+u(t)-f_2'(y(t))-f_1'(r^*)\le 0$, which implies that the integrand is nonpositive also
in this case, as claimed. In conclusion, the expression on the \rhs\ of \eqref{FabFour} is nonpositive.
At this point, we integrate \eqref{FabFour} over $\pier{(0,t)}$, where $t\in (0,T]$ is arbitrary. Since 
$(y_0-r^*)^+=0$ by assumption, we obtain that $(y-r^*)^+=0$ a.e. in $Q$, which implies
$y\le r^*$ a.e. in $Q$. Similarly, we obtain that $y\ge r_*$ a.e. in $Q$. 
\revis{With this, the validity of {\bf (GB)} is shown}.

\juerg{%
We conclude this \revis{example} with the remark that the above argumentation remains valid for every potential 
$f_1\in C^1(-1,1)\cap C^0([-1,1])$ which is convex on $[-1,1]$ and satisfies
\begin{equation*}
\lim_{r\searrow -1}f_1'(r)=-\infty, \quad \lim_{r\nearrow +1}f_1'(r)=+\infty,
\end{equation*}
where it is understood that $f_1$ is extended to the whole of $\erre$ by putting $f_1(r)=+\infty$
for $r\not\in [-1,1]$.}  

\vspace{1mm}
{\bf 2.} \,\,Next, we assume that $f_1\in C^3(\erre)$, which is satisfied for the classical potential \eqref{regpot}. 
In this case, $\VB\sigma\subset H^{2\sigma}(\Omega)$,
and it holds $H^{2\sigma}(\Omega)\subset L^\infty(\Omega)$ 
(and thus $y\in L^\infty(Q)$ with \eqref{bounds2} whenever $(\mu,y)={\cal S}(u)$ for some
$u\in{\cal U}_R$) if $\sigma>3/4$. 

\vspace{1mm}
{\bf 3.} The following result shows that the \revis{assumption} $\sigma>3/4$ 
\revis{(which ensures boundedness for $y$ whenever $f_1\in C^3(\erre)$ as shown in Example~2)
is not optimal if the nonlinearity satisfies an additional growth condition}. 
We remark that \revis{also} this condition is met by, e.g., the classical regular potential~\eqref{regpot}.

\Bprop 
Let $B=-\Delta$ with domain $H^2(\Omega)\cap H_0^1(\Omega)$, let $f\in C^3(\erre)$, and suppose
that the general assumptions {\bf (A1)}--{\bf (A5)} and {\bf (A7)} are
fulfilled. In addition,  assume that there is some $\widehat C_1>0$ such
that 
\begin{equation}\label{growth2}
|f'(s)|\,\le\,\widehat C_1\left(1+|s|^3\right) \quad\forall\,s\in\erre.
\end{equation}
Then the condition {\bf (GB)} holds true  whenever 
\,\,$\frac 9{20} < \sigma\le \frac 34$. 
\Eprop

\Bdim
We  show the result only for $\frac 9{20}<\sigma
<\frac 34$ (the case $\sigma=\frac 34$ can be treated in a 
similar way). We then have  
\begin{align}\label{findp}
\VB\sigma\subset H^{2\sigma}(\Omega)\subset L^p(\Omega)
\quad\mbox{with \,}-\frac 3p=2\sigma-\frac 32,\,\,\mbox{ \,i.e., \,}\,\,p=\frac 6{3-4\sigma}\,.
\end{align}
\oldgianni{We notice that \eqref{findp} holds true also in the case of Neumann boundary conditions.
However, we have to assume Dirichlet boundary conditions later on}.
From \eqref{Iseconda}, we infer that $\,B^{2\sigma}y=g-f'(y)$\, with \,$g:=\mu+u-\tau\,\dt y$,
where, owing to \eqref{bounds1} and \eqref{growth2},
\begin{equation}\label{SirToby}
\|g\|_{L^\infty(0,T;H)}\,+\,\|f'(y)\|_{L^\infty(0,T;L^{p/3}(\Omega))}\,\le\,C_1\,,
\end{equation}
with a global constant $C_1>0$. We now distinguish between the two cases $\,p/3\ge 2\,$ and
$\,p/3<2$, which, by virtue of \eqref{findp}, occur if $\,\sigma\ge 1/2\,$ and
$\,\sigma<1/2$, respectively.

Assume first that $\sigma\ge 1/2$. Then, by \eqref{SirToby}, $\,B^{2\sigma}y\in L^\infty(0,T;H)$,
whence
$$ y\in L^\infty(0,T;V_B^{2\sigma})\subset L^\infty(0,T;H^{4\sigma}(\Omega)) \subset L^\infty(Q),$$
since $\,4\sigma\ge 2$. Therefore, \eqref{bounds2} is valid.

Assume now that $\sigma<1/2$. Then, we only have $\,B^{2\sigma}y\in L^\infty(0,T;L^{p/3}(\Omega))$.
We now claim that the following implication is valid:
\begin{align}\label{claim}
&\mbox{If $\,\,v\in H\,\,$ and \,\,$B^s v\in L^q(\Omega)$\,\, with \,\,$s\in (0,1)$\,\, and
\,\,$q>\frac 3{2s}$, \,\,then \,\,$v\in L^\infty(\Omega)$}\nonumber\\
& \mbox{and \,\,$\|v\|_{\oldgianni\Linfty}\,\le\,C_2$, \,where\,\,$C_2$\,\, depends only on \,\,$s$, $q$ \,and\,
$\Omega$.}  
\end{align}  
To prove this claim, we note that $\lambda_1'>0$ \oldgianni{(see~\eqref{eigen})} in our situation, 
and thus we have 
$\,B^sw=0\,$ for $\,w\in V_B^\sigma\,$ if and only if $\,w=0$. Therefore, 
we must have $v=\tilde v_+-\tilde v_-$, where $\,\tilde v_\pm\in V_B^\sigma\,$ is the 
(unique) weak solution 
to the fractional Dirichlet problem $\,B^s \tilde v_\pm=(B^s v)_\pm$. 
At this point, \oldgianni{as we are dealing with Dirichlet boundary conditions},
we can apply the results of \cite[Thm.~4.1 and Sect.~2.1]{BoFiVa}, which 
imply that the estimate
\begin{equation}\label{claimy}
0\,\le\, \tilde v_\pm\,\le\,\kappa\,\left\|(B^sv)_\pm\right\|_q\,{\calB}_q(\phi) \quad \mbox{in }\,\Omega,
\end{equation}
holds true. Here, the constant $\kappa>0$ depends on $s$, $q$, and $\Omega$, $\phi\in C^0(\bar\Omega)$ 
is the first (positive)
eigenfunction of \,$B^s$\, (or, equivalently, of  $\,B$, i.e., we have $\phi=e_1'$), and $\,{\cal B}_q\,$ is a suitable
continuous function on $\,[0,+\infty)$\, depending on $\,q$. The claim thus holds true. 

We now choose $\,s=2\sigma$, so that $s\in (0,1)$, as well as $q=p/3$. Then we can apply \eqref{claim}
provided that $q>\frac 3{2s}$, \,\,i.e.,\,\, $2s>\frac 3q$, which, in view of \eqref{findp}, just means that
$\,\sigma>\frac 9{20}$. 
\Edim

\Brem
Observe that if  $B=-\Delta$ with zero Dirichlet boundary condition
and $\sigma>\frac 9{20}$, then the assumption {\bf (A2)} can only be fulfilled if 
$\lambda_1>0$. Indeed, if $\lambda_1=0$, then {\bf (A2)},(ii) necessitates that the constant
functions belong to $\VB\sigma\subset H^{2\sigma}(\Omega)$, which in turn requires that
$0<\sigma<1/4$.
\Erem

\vspace{2mm}
\revis{In the result stated below, 
which improves the stability estimate \eqref{stabu1} established in Theorem~\ref{Wellposedness},
we assume that the conditions {\bf (GB)} and {\bf(A8)} are satisfied
and account for Remarks~\ref{CtoSmap} and~\ref{Addreg}}.

\Bthm
\label{Stability}
Suppose that {\bf (A1)}--{\bf (A5)}, \revis{{\bf (A7)}--{\bf(A8)}} and {\bf (GB)} are satisfied. 
Then there is a constant $K_4>0$, which
depends only on the data of the state system and $R$, such that the following holds true:
whenever $u_i\in{\cal U}_R$, $i=1,2$, are given and $(\mu_i,y_i)={\cal S}(u_i)$, $i=1,2$, are
the associated solutions to the state system \eqref{Iprima}--\eqref{Icauchy}, then it holds,
for every $t\in (0,T]$,
\begin{align}
\label{stabu2}	
&\|\mu_1-\mu_2\|_{L^2(0,t;V_A^{2r})}\,+\,\|y_1-y_2\|_{H^1(0,t;H)\cap L^\infty(0,t;
V_B^\sigma)}\nonumber\\
&\le\,K_4\,\|u_1-u_2\|_{L^2(0,t;H)}\,. 
\end{align}
\Ethm

\Bdim
The functions $u:=u_1-u_2$, $y:=y_1-y_2$, $\mu:=\mu_1-\mu_2$, obviously satisfy the system
\begin{align}\label{diff1}
&\dt y+A^{2r}\mu=0 \quad\mbox{a.e. in }\,Q,\\
\label{diff2}
&\tau\,\dt y+B^{2\sigma}y +f'(y_1)-f'(y_2)=\mu+u\quad\mbox{a.e. in }\,Q,\\
&y(0)=0 \quad\mbox{a.e. in }\,\Omega.\label{diff2.5}
\end{align}
In the following, $C_i$, $i\in\enne$, denote constants that depend only on the data of
the state system and $R$.
We multiply \eqref{diff1} by $\mu$ and \eqref{diff2} by $\dt y$, add the resulting identities,
and integrate over $Q_t$, where $t\in (0,T]$ is arbitrary. Rearranging terms and applying
Young's inequality, we then obtain  
the inequality
\begin{align}
\label{diff3}
&\tau\int_0^t\!\|\dt y(s)\|^2\,ds\,+\int_0^t\|A^r\mu(s)\|^2\,ds\,+\,\frac 12\|B^\sigma y(t)\|^2
= \int_0^t\!\!\iO \dt y\bigl(u-(f'(y_1)-f'(y_2)\bigr)\nonumber
\\
&\le \frac\tau 2\int_0^t\!\|\dt y(s)\|^2\,ds\,+\,C_1\,\left(\|u\|_{L^2(Q_t\oldgianni)}^2
\,+\,\|f'(y_1)-f'(y_2)\|^2_{L^2(Q_t)}\right).
\end{align}
Now observe that $\,|f'(y_1)-f'(y_2)|\le K_3|y|\,$ a.e. in $Q$, by \eqref{bounds2}.
Hence, if we add the term $\,\int_{Q_t} y\,\dt y\,$ to both sides of \eqref{diff3} and apply
Young's inequality appropriately, then we readily infer from Gronwall's lemma the
estimate
\begin{align}
\label{stabu3}
&\|A^r\mu\|_{L^2(0,t;H)}\,+\,\|y\|_{H^1(0,t;H)\cap L^\infty(0,t;V_B^\sigma)}
\,\le\,C_2\,\|u\|_{L^2(0,t;H)}\,,
\end{align}
whence, by virtue of \eqref{diff1},
also
\begin{equation}\label{stabu4}
\|A^{2r}\mu\|_{L^2(0,t;H)}\,\le\,C_2\,\|u\|_{L^2(0,t;H)}\,.
\end{equation}
It remains to show the estimate
\begin{align}
\label{stabu5}
\|\mu\|_{L^2(0,t;V_A^{2r})}\,\le\,C_3\,\|u\|_{L^2(0,t;H)}\,.
\end{align}
According to \eqref{defnormaAr}, this follows directly from \eqref{stabu4} if
$\lambda_1>0$, while in the case $\lambda_1=0$ we have to estimate the mean value 
$\,\,{\rm mean}\,(\mu)$. Now, by {\bf (A2)}, 
the constant function $\,{\bf 1}(x)\equiv 1$\, belongs to $V_B^\sigma$. 
Moreover, we have in this case that $A^r{\bf 1}=0$,
and it follows from \eqref{diff1} that $\,{\rm mean} \,(\dt y)=0$, almost everywhere on $(0,T)$. 
We thus can integrate \eqref{diff2} over $\Omega$ to see
that we have almost everywhere in $(0,T)$ the estimate
\begin{align*}
\Bigl|\iO\mu(t)\Bigr|\,&\le\,\Bigl|\tau\iO\dt y(t)\,+\,(B^\sigma\pier{y(t)},B^\sigma{\bf 1})
-\iO u(t) +\iO(f'(y_1(t))-f'(y_2(t)))\Bigr|\nonumber\\[1mm]
&\le\,C_4\left(\|B^\sigma y(t)\|+\|u(t)\|+\|y(t)\|\right),
\end{align*}
and \eqref{stabu3} implies that
$$\|{\rm mean}\,(\mu)\|_{L^2(0,t)}\,\le\,C_5\,\|u\|_{L^2(0,t;H)},$$
whence \eqref{stabu5} follows.

\Edim

\section{Differentiability of the control-to-state mapping}
\label{FRECHET}
\setcounter{equation}{0}
In this section, we prove that the control-to-state mapping ${\cal S}:u\mapsto (\mu,y)$ is Fr\'echet
differentiable \oldgianni{from the space $\calX$ \pier{defined} in \eqref{defX} into a suitable Banach space~$\calY$}. 
To this end, we \revis{suppose} that 
the general assumptions {\bf (A1)}--{\bf (A5)}, {\bf (A7)} and {\bf (GB)} are satisfied.
\revis{By Remark~\ref{CtoSmap} the \revis{control}-to-state map $\cal S$ is well defined}.
We fix some $\bar u\in {\cal U}_R$ and let $( \bar \mu, \bar y)={\cal S}(\bar u)$. 
We then consider for an arbitrary $k\in{\cal X}$ the linearized system
\begin{align}
\label{lin1}
&\dt\xi+A^{2r}\eta=0 \quad\mbox{in $Q$},
\\[1mm]
\label{lin2}
&\tau\,\dt\xi+B^{2\sigma}\xi+f''(\bar y)\xi = \eta+k \quad\mbox{in $Q$},
\\[1mm]
\label{lin3}
&\xi(0)=0 \quad\mbox{in $\Omega$}. 
\end{align}
More precisely, we consider \revis{its weak version}
\begin{align}
\label{wlin1}
&(\dt\xi(t),v)\,+\left(A^r\eta(t),A^r v\right)\,=\,0\quad\mbox{for a.e. $\,t\in (0,T)\,$ and all $\,v\in V_A^r$},\\[2mm]
\label{wlin2}
&(\tau\,\dt \xi(t),v)\,+\left(B^\sigma \xi(t),B^\sigma v\right)\,+\,(f''(\bar y(t))\xi(t),v)\,=\,
(\eta(t)+k(t),v)\nonumber\\[1mm]
&\quad \mbox{for a.e. \,$t\in (0,T)$\, and all \,$v\in V_B^\sigma$},\\[2mm]
\label{wlin3}
&\xi(0)=0.
\end{align}
If this system admits a unique solution $(\eta,\xi)$, and if the Fr\'echet derivative
$\,D{\cal S}(\bar u)\,$ of $\,{\cal S}\,$ at \,$\bar u$\, exists, then we should have 
that $\,D{\cal S}(\bar u)(k)=(\eta,\xi)$. Observe that $\,\bar y\,$ enjoys the 
regularity \eqref{regy}, and the global
bounds \eqref{bounds1} and \eqref{bounds2} are satisfied for $y=\bar y$. We have the following result.

\Bthm
\label{Linearized}
Under the given assumptions, the linearized system \eqref{wlin1}--\eqref{wlin3} admits for 
every $\bar u\in\uad$ and every
{$k\in {\cal X}$} a unique solution $(\eta,\xi)$ such that 
\begin{align}
\label{reglin}
\eta\in L^2(0,T;V_A^r), \quad \xi\in H^1(0,T;H)\cap L^\infty(0,T;V_B^\sigma).
\end{align}
Moreover, there is a constant $K_5>0$, 
which depends only on the data of the state system and $R>0$, such that
\begin{align}
\label{stabulin}
\oldgianni{%
  \norma\eta_{\L2{\VA r}}
  + \norma\xi_{\H1H\cap\L\infty{\VB\sigma}}
  \leq K_5 \, \norma k_{\LQ\infty} \,.
}%
\end{align}
\Ethm

\Bdim
We prove the assertion in a number of \revis{steps}.

\step
Step 1. Discretization

We fix an integer $N>1$,  set $h_N:=T/N$ and $\,t_N^n:=n\,h_N$, $n=0,\ldots,N$,
 and notice that by virtue of the global bound
\eqref{bounds2} the linear operators
\begin{equation}
\label{defP}
P_N^n:H\to H; \quad v\mapsto P_N^n\,v:=f''(\bar y(\cdot,t_N^n)v,
\end{equation}
are continuous, where with $\widehat C:=K_3$ it holds
\begin{equation*}
\|P_N^n\|_{{\cal L}(H,H)}\,\le\,\widehat C\quad\,\forall\, N\in \enne, \quad 0\le n\le N.
\end{equation*}
The discrete problem then consists in finding
two $\,(N+1)$-tuples\, $(\xi^0_N,\dots,\xi^N_N)$ and $(\eta^0_N,\linebreak\dots,\eta^N_N)$
satisfying
\Beq
  \xi^0_N = \eta^0_N=0\,, \quad
  (\xi^1_N,\dots,\xi^N_N) \in (\VB{2\sigma})^N,
  \quad
  (\eta^1_N,\dots,\eta^N_N) \in (\VA{2r})^N,
  \label{regdiscr}
\Eeq
and 
\begin{align}
  & \dhyn + \munp + A^{2r} \munp
  \,=\, \mun,
  \label{primad}
  \\[1mm]
  & \tau \, \dhyn
  + \left(\widehat C I + B^{2\sigma} + P_N^{n+1}\right)(\ynp)
   \,=\, \widehat C\, \yn + \munp + k^{n+1}_N,
  \label{secondad}
\end{align}
for $n=0,1,\dots,N-1$, where $I:H\to H$ is the identity and
\Beq
  k^n_N := k(nh_N)
  \quad \hbox{for $n=0,1,\dots,N$}.
  \label{defun}
\Eeq
\pier{In view of \eqref{defX}, note that $k$ is continuous from $[0,T]$ to $H$, so
that the above definition is meaningful.
The problem \eqref{regdiscr}--\eqref{secondad}}
can be solved inductively for $n=0,\dots,N-1$ in the following way:
let $(\mun,\yn)$ be given in $\VA{2r}\times \pier{\VB{2\sigma}}$.
We first rewrite the above equations in the form
\begin{align}
  & h_N \left(I+A^{2r}\right) \munp + \ynp
  \,=\, \yn + h_N \mun, 
  \label{primadbis}
  \\[1mm]
  & \left(\oldgianni{(\widehat C+(\tau/h_N))I} + B^{2\sigma} + P_N^{n+1}\right)(\ynp)
   \,=\, (\widehat C+(\tau/h_N)) \yn + \munp + k^{n+1}_N.
  \label{secondadbis}
\end{align}
Next, we observe that the operator $\,\calA:=\widehat C\,I+P_N^{n+1}:H\to H\,$
is monotone and continuous.
On the other hand, the unbounded operator $\,B^{2\sigma}\,$~is monotone in $H$,
and $I+ B^{2\sigma}:V_B^{2\sigma}\to H$ is surjective, whence it follows that
 $B^{2\sigma}$ is maximal monotone.
Therefore, the sum $\,\calA+B^{2\sigma}\,$ is also maximal monotone
(see, e.g., \cite[Cor.~2.1 p.~35]{Barbu}).
It follows that $\,(\tau/h_{N})I+\calA+B^{2\sigma}$, 
i.e., the operator that acts on $\ynp$ in~\eqref{secondad},
is \pier{linear,} surjective and one-to-one from $\VB{2\sigma}$ onto~$H$.
Therefore, \eqref{secondad}~can be rewritten in the equivalent form
\Beq
  \ynp
  = \left(L_NI + B^{2\sigma} + P_N^{n+1}\right)^{-1} \bigl(L_N \yn + \munp + k^{n+1}_N),
  \label{secondadter}
\Eeq
where, for brevity, we have set\, $L_N:=\widehat C+(\tau/h_N)$.
By accounting for \eqref{primadbis}, we conclude that problem \accorpa{primad}{secondad}
is equivalent to the system obtained by coupling \eqref{secondadter} with
the equation
\Beq
  h_N\left(I+A^{2r}\right) \munp
  + \left(L_N I + B^{2\sigma} + P_N^{n+1}\right)^{-1} \bigl(L_N \yn + \munp + k^{n+1}_N)
  \,=\, \yn + h_N \mun 
  \non
\Eeq
or
\Bsist
  && h_N\left(I+A^{2r}\right) \munp
  + \left(L_N I + B^{2\sigma} + P_N^{n+1}\right)^{-1} \munp
  \non
  \\
  && =\, \yn + h_N \mun 
  - \left(L_N I + B^{2\sigma} + P_N^{n+1}\right)^{-1} \bigl(L_N \yn + k^{n+1}_N) .
  \label{secondadquater}
\Esist
By arguing as before, we see that the operator acting on $\munp$ 
on the \lhs\ of \eqref{secondadquater} is surjective and one-to-one from $\VA{2r}$ onto~$H$,
so that the equation can be uniquely solved for~$\munp$ in~$\VA{2r}$.
Inserting the solution in \eqref{secondadter}, 
we directly find that $\ynp\in\VB{2\sigma}$.

\medskip

Now that the discrete problem is solved, we can start \revis{to perform the a priori estimates.}
In the following, the (possibly different) values of the constants termed $C_i$, $i\in \enne$, 
are independent of the parameters $h_N=T/N$ and~$n\in\enne$. 
Also,  in order to avoid an overloaded notation, we omit the index $N$ in the expressions $\xi_N^n$ and $\eta_N^n$, writing
it only at the end of each estimate.

Moreover, we also express the bounds we find in terms of the interpolants.
According to the notation introduced in Section 2, and recalling that $\xi_N^0=\eta_N^0=0$,
we remark at once that the discrete problem also reads
\begin{align}
  & \widehat\xi_{h_N} \in \W{1,\infty}{\VB\sigma} , \quad
  \underyh,  \overyh \in \L\infty{\VB{2\sigma}},
  \qquad
  \label{regyh}
  \\[1mm]
	  & \undermuh \,, \overmuh \in \L\infty{\VA{2r}},
  \label{regmuh}
  \\[1mm]
  & \dt\yh
  + \overmuh
  + A^{2r} \overmuh
  \,=\, \undermuh
  \quad \aet,
  \label{primah}
  \\[1mm]
  & \tau \, \dt\yh
  + \left(\widehat C I + B^{2\sigma} + P_{N}\right)(\overyh)
  \, = \,\widehat C \underyh + \overmuh + \overuh
  \quad \aet ,
  \label{secondah}
  \\[1mm]
  & \yh(0) = 0 \,,
  \label{cauchyh}
\end{align}
where it is understood that 
\begin{equation}
\label{defPN}
\left(P_N \overyh\right)(\cdot,t)=P_N^{n+1}\xi_N^{n+1} \quad \mbox{for a.e. 
$\,t\in (t_N^n,t_N^{n+1}), \quad 0\le n\le N-1$.}
\end{equation}

\step
Step 2. First a priori estimate

We test \eqref{primad} and~\eqref{secondad} 
(by~taking the scalar product in~$H$)
by~$h_N \eta^{n+1}$ and $\xi^{n+1}-\xi^n$, respectively,
and add the resulting identities.
Noting an obvious cancellation, we obtain the equation
\begin{align}\label{j0}
  & h_N\, (\eta^{n+1}-\eta^n,\eta^{n+1})
  + h_N\, (A^{2r} \eta^{n+1},\eta^{n+1})
    + \, \tau \,h_N \, \left\|\mbox{$\frac{\xi^{n+1}-\xi^n}{h_N}$}\right\|^2
	\non\\[1mm]	
  &+ (B^{2\sigma}\xi^{n+1},\xi^{n+1}-\xi^n)
    + \bigl( (\widehat C I + P_N^{n+1}) (\xi^{n+1}),\xi^{n+1}-\xi^n \bigr)
  \non
  \\[1mm]
  & =\, \widehat C (\xi^n,\xi^{n+1}-\xi^n)
  + (k^{n+1},\xi^{n+1}-\xi^n).
  \end{align}
Now, we observe that 
\begin{equation}\label{j1}
\left(\widehat C\,\xi^{n+1},\xi^{n+1}-\xi^n\right)\,=\,\frac {\widehat C}2\,\|\xi^{n+1}\|^2-\frac{\widehat C}2\,\|\xi^n\|^2
+\frac{\widehat C}2\,\|\xi^{n+1}-\xi^n\|^2\,.
\end{equation}
Moreover, by Young's inequality it holds  that
\begin{align}\label{j2}
&\left(\widehat C\,\xi^n \,+\,k^{n+1}\,-\,P_N^{n+1}\xi^{n+1},\xi^{n+1}-\xi^n\right)\non\\
&\le\,\frac \tau 2\,h_N\,\left\|\mbox{$\frac{\xi^{n+1}-\xi^n}{h_N}$}\right\|^2\,+\,
\frac 1{2\tau}\,h_N\left\|\widehat C\,\xi^n\,+\,k^{n+1}\,-\,P_N^{n+1}\xi^{n+1}\right\|^2\non\\
&\le\,\frac \tau 2\,h_N\,\left\|\mbox{$\frac{\xi^{n+1}-\xi^n}{h_N}$}\right\|^2\,+\,
C_1\,h_N\left(\|\xi^n\|^2\,+\,\|k\|_{L^\infty(Q)}^2\,+\,\|\xi^{n+1}\|^2\right),
\end{align}
where $C_1$ depends only on $\tau$ and $\widehat C$. 
Combining \eqref{j0}--\eqref{j2}, 
we deduce that
\begin{align}\label{j3}
& \frac {h_N}2 \, \|\eta^{n+1}\|
  - \frac {h_N}2 \,\|\eta^n\|^2 + \frac{h_N}2\,\|\eta^{n+1}-\eta^n\|^2
    + h_N\,\|A^r \eta^{n+1}\|^2
  \non
  \\
  & \quad {}
  + \frac \tau 2\, h_N\left\|\mbox{$\frac{\xi^{n+1}-\xi^n}{h_N}$}\right\|^2
  + \frac 12 \, \norma{B^\sigma \xi^{n+1}}^2
  + \frac 12 \, \norma{B^\sigma(\xi^{n+1}-\xi^n)}^2
  - \frac 12 \, \norma{B^\sigma\xi^n}^2
  \non
  \\
  & \quad {}
  + \frac {\widehat C} 2 \, \|\xi^{n+1}\|^2 - \frac{\widehat C}2\,\|\xi^n\|^2
  +\frac{\widehat C}2\,\|\xi^{n+1}-\xi^n\|^2 \non\\
  &\le \,C_1\,h_N\left(\|\xi^n\|^2\,+\,\|k\|_{L^\infty(Q)}^2\,+\,\|\xi^{n+1}\|^2\right).
  \end{align}
Then, we sum up for $n=0,\dots,\ell-1$ with $\ell\leq N$, obtaining the inequality
\begin{align}
  & \frac {h_N}2 \, \|\eta^\ell\|^2 
  + \oldgianni{\frac 12 \,\somma n0{\ell-1} {h_N} \, \|\eta^{n+1}-\eta^n\|^2}
  + \somma n0{\ell-1} h_N\,\|A^r \eta^{n+1}\|^2
  \non
  \\
  & \quad {}
  + \frac \tau 2 \somma n0{\ell-1} h_N \left\|\mbox{$\frac {\xi^{n+1}-\xi^n}{h_N}$}\right\|^2
  + \frac 12 \, \norma{B^\sigma \xi^\ell}^2
  + \frac 12 \, \somma n0{\ell-1} \norma{B^\sigma(\xi^{n+1}-\xi^n)}^2
  \non
  \\
  & \quad {}
  + \Bigl(\frac {\widehat C}2-C_1\,h_N\Bigr)\,\|\xi^\ell\|^2
  + \frac{\widehat C}2 \,\sum_{n=0}^{\ell-1}\|\xi^{n+1}-\xi^n\|^2
  \non
  \\[1mm]
  & \leq \,C_1\,\ell\,h_N\,\|k\|_{L^\infty(Q)}^2 \,+\,2\,C_1\,\sum_{n=0}^{\ell-1}
  h_N\,\|\xi^{n}\|^2\,.
  \label{perprimastima}
\end{align}
At this point, we fix any $N_0\in\enne$ such that \,$N_0\ge\oldgianni{4\,C_1\,T/\widehat C}$. With this choice,
we have for any integer $\,N\ge N_0$\, that $\,\frac{\widehat C}2-C_1\,h_N\ge\frac {\widehat C}4$. Since also
$\ell\,h_N\le T$, we conclude from the discrete Gronwall lemma that for any such $N\in\enne$ it holds the
bound
\begin{align}
\label{j4}
& \pier{h_N} \, \|\eta^\ell\|^2
  + \pier{\somma n0{\ell-1} {h_N} \, \|\eta^{n+1}-\eta^n\|^2}
  + \somma n0{\ell-1} h_N\, \|A^r\eta^{n+1}\|^2
  + \somma n0{\ell-1} h_N \,\left\|\mbox{$\frac{\xi^{n+1}-\xi^n}{h_N}$}\right\|^2
  \non
  \\
  & \quad {}
  + \, \|B^\sigma \xi^\ell\|^2
  + \,\somma n0{\ell-1} \, \norma{B^\sigma(\xi^{n+1}-\xi^n)}^2 \,+\,\|\xi^\ell\|^2\,+\,\sum_{n=0}^{\ell-1}
  \|\xi^{n+1}-\xi^n\|^2
    \non
  \\
& \leq\, C_2\,\|k\|_{L^\infty(Q)}^2\,\le\,C_3 \,.
\end{align}
\pier{Since this holds for $\ell=0,\dots,N$, 
%(where we put $\xi^{N+1}_N=\xi^N_N$ and $\eta_N^{N+1}=\eta_N^N$), 
we obtain} in terms of the interpolants, by neglecting the first contribution
and recalling that $\mu^0=0$ \oldgianni{and the definition~\eqref{defBs} of the norm in~$\VB\sigma$}, that
\begin{align}
  & \norma{\overmuh-\undermuh}_{\L2H}
  + \norma{A^r \overmuh}_{\L2H}
  + \norma{A^r \undermuh}_{\L2H}
  \pier{{}+ \norma{\dt\yh}_{\L2H}}  \non
  \\[2pt] 
  & \quad {}
  + \oldgianni{\norma\underyh_{\L\infty{\VB\sigma}}
  + \norma\overyh_{\L\infty{\VB\sigma}}}
%  + \norma\yh_{\L\infty{\VB\sigma}}}
%  \non
%  \\
%  & \quad {}
  \pier{{}+ h_N^{-1/2}
  % \left( \norma{B^\sigma(\overyh-\underyh)}_{\L2H}
  \norma{\overyh-\underyh}_{\L\infty{\VB\sigma}}}
  % \right)
  \non
  \\
  & \le\,C_4\,\|k\|_{L^\infty(Q)}\,\le\,C_5\,.
    \label{primastima}
\end{align}

\step
Step 3. Second a priori estimate

Let $N\ge N_0$. We want to improve the estimate for $A^r\overmuh$ given by~\eqref{primastima}
and show that
\Beq
  \norma\overmuh_{\L2{\VA r}} 
  + \norma\undermuh_{\L2{\VA r}} \leq C_6 \|k\|_{L^\infty(Q)}\leq C_7\,.
  \label{terzastima}
\Eeq
By recalling \eqref{defnormaAr},
we see that there is nothing to prove if $\lambda_1>0$. Assume now that
$0=\lambda_1<\lambda_2$. We then have to estimate the mean value of~$\overmuh$.
To this end, we recall that \,$e_1\,$ is a constant and belongs to $\,V_B^\sigma$.
Thus, the function ${\bf 1}(x)\equiv 1$ also belongs to $\,V_B^\sigma$. Integrating
the \oldgianni{equation}~\eqref{secondah} over~$\Omega$, we therefore obtain almost everywhere on $(0,T)$ the 
identity
\begin{align}
\iO \overmuh \,=\,\iO\Bigl(-\overuh + \widehat C\left(\overyh-\underyh\right)+ P_N(\overyh)
+\,\tau\,\partial_t\yh\Bigr) + \left(B^\sigma \overyh,B^\sigma {\bf 1}\right).
\end{align}
Applying the Cauchy--Schwarz inequality to the expressions on the right-hand side, we
readily conclude from \eqref{primastima} the bound
\begin{align}
\left\|\mbox{mean}(\overmuh)\right\|^2_{L^2(0,T)} \,&\le\,C_8 \Bigl(\|k\|_{L^\infty(Q)}^2
\,+\,\|\overyh\|_{L^2(Q)}^2\,+\,\|\underyh\|_{L^2(Q)}^2 \,+\,\left\|B^\sigma \overyh\right\|^2_{L^2(Q)}
\non\\
&\hspace*{12mm} + \,\bigl\|\partial_t\yh\bigr\|^2_{L^2(Q)}\Bigr)\non\\
&\le\,C_9\,\|k\|^2_{L^\infty(Q)}\,\le\,C_{10}\,,
\end{align} 
and the claim \eqref{terzastima} is proved as far as $\,\overmuh\,$ is concerned. But as $\,\|\overmuh-\undermuh\|_{L^2(Q)}\,$
is by \eqref{primastima} bounded, and since $\,A^r\eta^0_N=A^r0=0$, it also holds true for $\,\undermuh$.

\step
Step 4. Existence

Combining  the estimates \eqref{primastima} and \eqref{terzastima}, \pier{recalling 
\eqref{diffLdue},}
and using standard weak and weak-star compactness results, we see that there are functions
$\xi$ and $\eta$ such that, at least for suitable subsequences which are again indexed by $N$,
\begin{align}
\label{convyh}
  & \overyh \to \xi \,, \quad
  \underyh \to \xi\,,
  \quad
  \yh \to \xi,
  \quad \hbox{all weakly star in $\L\infty{\VB\sigma}$},
   \\[1mm]
 \label{convdtyh}
  & \dt\yh \to \dt\xi
  \quad \hbox{weakly in $\L2H$},
  \\[1mm]
  & \overmuh \to \eta
  \quad \hbox{weakly in $\L2{\VA r}$},
  \label{convmuh}
\end{align}
\oldgianni{as $N\to\infty$}. 
Moreover, owing to the compact embedding $\VB\sigma\subset H$ (see~\eqref{compembB})
and to well-known strong compactness results (see, e.g., \cite[Sect.~8, Cor.~4]{Simon}), we obtain
from \accorpa{convyh}{convdtyh}  that
\Beq
  \yh \to \xi
  \quad \hbox{strongly in $C([0,T];H)$},
  \label{strongyh}
\Eeq
whence it follows that $\xi(0)=0$ and, using \eqref{diffLdue},
\begin{equation}
\label{strongoveryh}
\overyh\to\xi, \quad \underyh\to\xi, \quad\mbox{both strongly in $\,L^2(0,T;H)$.}
\end{equation}

Next, we prove that
\Beq
  \undermuh \to \eta
  \quad \hbox{weakly in $\L2{\VA r}$}.
  \label{convundermuh}
\Eeq
By \eqref{primastima} and \eqref{convmuh}, it suffices to check that
\Beq
  {}_{\L2{\VA{-r}}}\<v , \overmuh-\undermuh>_{\L2{\VA r}} \to 0
  \quad \hbox{as $N\to\infty$},
  \label{perconv}
\Eeq
for every $v$ belonging to a dense subspace \oldgianni{$\calV$} of $\L2{\VA{-r}}$,
where we can take $\oldgianni\calV=\pier{C^1_c}(0,T;H)$ 
since $H$ is dense in~$\VA{-r}$ (see~\eqref{compembAneg}).
So, we fix $v\in \pier{C^1_c}(0,T;H)$ and choose $\delta>0$ such that 
$v(t)=0$ for $t\in[0,T]\setminus(\delta,T-\delta)$.
If $h_N\in(0,\delta/2)$, then we have
\begin{align}
  & |{}_{\L2{\VA{-r}}}\<v , \overmuh-\undermuh>_{\L2{\VA r}}|
  = \Bigl| \int_{h_N}^T (\overmuh-\undermuh)(t) \, v(t) \, dt \Bigr|
  \non
  \\
  & = \Bigl| \int_{h_N}^T \bigl( \overmuh(t) - \overmuh(t-h_N) \bigr) \, v(t) \, dt \Bigr|
  \non
  \\
  & = \Bigl| \int_{h_N}^T \overmuh(t) \, v(t) \, dt
  - \int_0^{T-h_N} \overmuh(t) \, v(t+h_N) \, dt \Bigr|
  \non
  \\
  & = \Bigl| \int_{h_N}^{T-h_N} \overmuh(t) \, (v(t)-v(t+h_N)) \, dt \Bigr|
  \leq \pier{T^{1/2}}\, \norma\overmuh_{\L2H}\, \norma{v'}_{\L\infty H} \, \pier{ h_N }\,,
  \non
\end{align}
and \eqref{perconv} follows.

We now show that
\begin{equation}
\label{conny}
P_N(\overyh)\to f''(\bar y)\xi \quad\mbox{strongly in \,\oldgianni{$L^1(Q)$} \,as\, $N\to\infty$.}
\end{equation}
Indeed, employing the global bounds \eqref{bounds2},
we have, for almost every $(x,t)\in \Omega\times (t_N^{n-1},t_N^n)$, where $1\le n\le N$,
\begin{align}
\label{yeah!}
&\left|P_N(\overyh)(x,t)-f''(\bar y(x,t))\xi(x,t)\right|=\,\left|f''(\bar y(x,t_N^n))\xi_N^n(x)-f''(\bar y(x,t))\xi(x,t)\right|\nonumber\\[1mm]
&\le\,\left|f''(\bar y(x,t_N^n))-f''(\bar y(x,t))\right||\xi(x,t)|+\left|f''(\bar y(x,t_N^n)\right|\left|\xi_N^n(x)-\xi(x,t)\right|\nonumber\\[1mm]
&\le\,\pier{\widehat C\left|\overyh(x,t)-\xi(x,t)\right|
\,+\,{}} \widehat C\,|\xi(x,t)|\left|\bar y(x,t_N^n)-\bar y(x,t)\right|\nonumber\\[1mm]
&\le\,\pier{\widehat C\left|\overyh(x,t)-\xi(x,t)\right|
\,+\,{}}\widehat C\,|\xi(x,t)|\int_{t_N^{n-1}}^{t_N^n}\left|\dt \bar y(x,s)\right|ds 
\nonumber\\[1mm]
&\le\,
\oldgianni{\widehat C\left|\overyh(x,t)-\xi(x,t)\right|
\,+\,\widehat C\,h_N^{1/2}\,|\xi(x,t)|\Bigl(\int_{t_N^{n-1}}^{t_N^n}|\dt\bar y(x,s)|^2\,ds\Bigr)^{1/2}}
\,. 
\end{align}
The claim \eqref{conny} then follows 
\oldgianni{from \eqref{strongoveryh}
and a simple calculation on the last term by recalling that $\xi$ and $\dt \bar y$ belong to~$L^2(Q)$}.

Therefore, we can pass to the limit as $N\to\infty$ in the weak time-integrated versions of \eqref{primah} and \eqref{secondah}  
(written with \oldgianni{bounded} time-dependent test functions) to conclude that the pair $(\eta,\xi)$
 solves the variational equations \eqref{wlin1} and~\eqref{wlin2}.
Since also $\xi(0)=0$, the existence part of the assertion is shown. Moreover, the continuity estimate \eqref{stabulin}
is a direct consequence of \eqref{primastima}, \eqref{terzastima} and the semicontinuity of norms.

\step
Step 5. Uniqueness

To show uniqueness, suppose that the system 
\eqref{wlin1}--\eqref{wlin3} has two solutions $(\eta_i,\xi_i)$, $i=1,2$, with the regularity \eqref{reglin}. Then the pair $(\eta,\xi)$ with $\eta=\eta_1-\eta_2$,
$\xi=\xi_1-\xi_2$,  solves the system \eqref{wlin1}--\eqref{wlin3}, where in this case $k\equiv 0$. We then test 
\eqref{wlin1} by $\eta$ and \eqref{wlin2} by $\dt \xi$ and add the resulting equations to arrive at the identity
\begin{align}
\label{uni}
\int_0^t\|A^r\eta(s)\|^2\,ds\,+\,\frac \tau 2\,\|B^\sigma \xi(t)\|^2 \,+\,\tau\int_0^t\!\!\iO|\dt\xi|^2\,=\,-\int_0^t\!\!\iO f''(\bar y)\,\xi\, \dt\xi\,,
\end{align}  
which is valid for every $t\in [0,T]$. Now we add the term $\,\,\int_0^t\!\!\iO \xi\,\dt\xi\,\,$ to both sides of \eqref{uni} and apply
Young's inequality appropriately to the resulting right-hand side. It then follows from Gronwall's lemma that 
$A^r\eta\,=\,\xi\,=\,0$. But then, by virtue of \eqref{wlin2}, also $\eta=0$.    
This concludes the proof of \oldgianni{Theorem~\ref{Linearized}}.
\Edim

After these preparations, the road is paved for proving the Fr\'echet differentiability of the
control-to-state operator ${\cal S}$. We need, however, yet another assumption.

\vspace{2mm}\noindent
{\bf \revis{(A9)}} \,\,\,$V_B^\sigma$ is continuously embedded in $L^4(\Omega)$.

\vspace{2mm}  \noindent
Observe that this condition is fulfilled if, e.g., $B=-\Delta$ with zero
Dirichlet \oldgianni{or Neumann boundary} conditions and $\sigma\ge 3/8$. 
Indeed, by virtue of \eqref{findp},
we have in this case  
$V_B^\sigma\subset H^{2\sigma}(\Omega)\subset L^4(\Omega)$ 
if \oldgianni{$-\frac34\leq2\sigma-\frac32$, i.e., if $\sigma\geq3/8$}.  

Recalling the statement of Theorem~\ref{Linearized}, we show the following result,
\revis{which also assumes {\bf(A8)} and accounts for Remarks~\ref{CtoSmap} and~\ref{Addreg}}.

\Bthm
\label{Differentiability}
Suppose that the assumptions {\bf (A1)}--{\bf (A5)}, \revis{{\bf (A7)}--{\bf (A9)}},
and {\bf (GB)} are fulfilled. 
 Then the control-to-state 
operator $\,{\cal S}: u\mapsto {\cal S}(u)=(\mu,y)\,$ is Fr\'echet 
differentiable in ${\cal U}_R$ when viewed as a mapping 
between the spaces $\,{\cal X}=H^1(0,T;H)\cap L^\infty(Q)\,$ 
and $\,{\cal Y}:=L^2(0,T;V_A^r)\times\left(H^1(0,T;H)\cap L^\infty(0,T;V_B^\sigma)\right).$
Moreover, whenever $\,\bar u\in {\cal U}_R\,$ with $(\bar\mu,\bar y)={\cal S}(\bar u)$ is given,
then the Fr\'echet derivative $D{\cal S}(\bar u)\in {\cal L}({\cal X},{\cal Y})$ \,of\, 
${\cal S}\,$ at \,$\bar u\,$ is \pier{specifed} by the identity $\,D{\cal S}(\bar u)(k)=(\eta,\xi)$,
where $(\eta,\xi)$ is the unique solution to the weak \pier{formulation \eqref{wlin1}--\eqref{wlin3} of the linearized system.}
%, which was established in \oldgianni{Theorem~\ref{Linearized}}. 
\Ethm
\Bdim
Since ${\cal U}_R$ is open, there is some $\,\Lambda>0\,$ such that $\bar u+k\in {\cal U}_R$
whenever $\,k\in{\cal X}\,$ and $\,\|k\|_{{\cal X}}\,\le\,\Lambda$. In the following, we
consider only such perturbations $\,k$, for which we define the quantities
\begin{align*}
&(\muk,\yk):={\cal S}(\bu+k), \quad \rk:=\muk-\bm-\etak, \quad \zk:=\yk-\by-\xik, 
\end{align*}
where $(\etak,\xik)=(\eta,\xi)$ denotes the unique solution to the system 
\eqref{wlin1}--\eqref{wlin3}. Obviously, we have \,$\oldgianni\rk\in L^2(0,T;V_A^r)\,$ and
\,$\oldgianni\zk\in \pier{ H^1(0,T;H) \cap {}}L^\infty(0,T;V_B^\sigma).$ Moreover, 
\pier{it turns out~that}
\begin{align}
\label{susi1}
&(\dt\zk(t),v)+\left(A^r \rk(t),A^r v\right)=0 \quad\mbox{for a.e. \,$t\in (0,T)$\, and all $\,v\in 
V_A^r$,}\\[2mm]
\label{susi2}
&\tau\,(\dt\zk(t),v)+\left(B^\sigma \zk(t),B^\sigma v \right)\,+\,
(f'(\yk(t))-f'(\by(t))-f''(\by(t))\xik(t),v)\nonumber\\[1mm]
&=\,(\rk(t),v)  \quad\mbox{for a.e. \,$t\in (0,T)$\, and all $\,v\in 
V_B^\sigma$,}\\[2mm]
\label{susi3}
&\zk(0)=0.
\end{align}
In addition, by Taylor's theorem and \eqref{bounds2}, we have almost everywhere in $Q$ that
\begin{align}\label{Taylor}
&\left|f'(\yk)-f'(\by)-f''(\by)\xik\right|\,\le\,C_1\left(|\zk|\,+\,|\yk-\by|^2\right),
\end{align}
where, here and in the remainder of the proof, the constants $C_i>0$, $i\in\enne$, depend
only on the data of the problem and $R$, but not on the special choice of $\,k\in{\cal X}\,$
with \,$\|k\|_{{\cal X}}\,\le\,\Lambda$. \pier{Using} \eqref{stabu2} in \oldgianni{Theorem~\ref{Stability}} 
and the continuity of the embedding $V_B^\sigma\subset L^4(\Omega)$, \pier{we infer} that\juerg{, for any $t\in (0,T]$,}
\begin{equation}
\label{susi4}
\|\yk-\by\|_{L^\infty(0,\oldgianni t;L^4(\Omega))}\,\le\,C_2\,\|k\|_{L^2(0,t;H)}\,.
\end{equation}

Now recall that by \eqref{stabulin} the mapping $\,\,k\mapsto (\etak,\xik)\,\,$ is 
continuous from \juerg{${\cal X}$ into ${\cal Y}$.} 
According to the notion of Fr\'echet differentiability, it therefore
suffices to construct an increasing function $\,Z:(0,\Lambda)\to (0,+\infty)\,$ such that
$\,\lim_{\lambda\searrow0} \frac {Z(\lambda)}{\lambda^2}\oldgianni{=0}$\, and
\begin{align}
\label{Frechet}
\|\rk\|^2_{L^2(0,T;V_A^r)}\,+\,\|\zk\|^2_{H^1(0,T;H)\cap L^\infty(0,T;V_B^\sigma)}
\,\le\,Z\left(\|k\|_{\L2H}\right)\,.
\end{align}

At this point, we test \eqref{susi1} by \,$\rk(t)$, \eqref{susi2} by \,$\dt \zk(t)$, add the
resulting equations, and integrate over $Q_t$, where $t\in (0,T]$. In addition, we add the 
term $\,\,\int_0^t\!\!\iO \zk\,\dt\zk \,\,$ to both sides of the result. Invoking
\eqref{Taylor}, we then obtain the
inequality
\begin{align}
\label{susi5}
&\frac 12\left(\|\zk(t)\|^2\,+\,\|B^\sigma \zk(t)\|^2\right)\,+\tau\int_0^t\!\!\iO|\dt\zk|^2
\,+\int_0^t\|A^r \rk(s)\|^2\,ds
\nonumber\\
&\le\,C_3\int_0^t\!\!\iO|\zk|\,|\dt\zk|\,+\,C_4\int_0^t\!\!\iO|\dt\zk|\left|\yk-\by\right|^2
\,=:\,I_1+I_2,
\end{align} 
with obvious notation. Now, by Young's inequality,
$$I_1\,\le\,\frac \tau 4\int_0^t\!\!\iO|\dt\zk|^2\,+\,C_5\int_0^t\!\!\iO|\zk|^2\,,$$
while, by also using H\"older's inequality and \eqref{susi4},
\begin{align*}
&I_2\,\le\,C_4\int_0^t\|\dt\zk(s)\|\,\|\yk(s)-\by(s)\|_{\oldgianni{\Lx4}}^2\,ds\,\le\,\frac \tau 4
\int_0^t\!\!\iO|\dt\zk|^2\,+\,C_6\,\|k\|^4_{\oldgianni{\L2H}}\,.
\end{align*}
Employing Gronwall's lemma, we thus conclude from \eqref{susi5} the estimate
\begin{align}
\label{susi6}
\|\zk\|^2_{\oldgianni{\H1H\cap\L\infty{\VB\sigma}}} \,+\,\|A^r\rk\|^2_{\oldgianni{\L2H}}
\,\le\,C_7\,\|k\|_{\oldgianni{\L2H}}^4\,.
\end{align}

At this point, we have to distinguish between two cases. Assume first that $\lambda_1>0$.
In this case, we have $\,\|\rk\|_{L^2(0,T;V_A^r)}\,\le\,C_8\,\|A^r\rk\|_{L^2(0,T;H)}$,
and thus \eqref{Frechet} follows from \eqref{susi6} with $\,Z(\lambda)=\oldgianni{(1+C_8)C_7\,\lambda^4}$.

Assume now that $\lambda_1=0$. In this case, we need to estimate the mean value of~$\rk$. 
To this end, we observe that {\bf (A1)} implies that for $\lambda_1=0$ we have 
${\bf 1}\in V_A^r\cap V_B^\sigma$ and $A^r{\bf 1}=0$. 
From this it immediately follows that \,$\mbox{mean}\,(\dt\zk(t))=0\,$ for almost every $t\in (0,T)$.  
Thus, inserting \,$v={\bf 1}\in V_B^\sigma\,$ in \oldgianni{\eqref{susi2}}
and applying the Cauchy--Schwarz inequality and the inequality~\eqref{Taylor}, 
we find that \aat\ it holds that
\begin{align*}
\Bigl|\iO\rk(t)\Bigr|\,&\le\,\iO |B^\sigma \oldgianni\zk(t)|\,|B^\sigma{\bf 1}|\,+\,
\oldgianni{C_1}\iO\left(|\oldgianni\zk(t)|\,+\,|\yk(t)-\by(t)|^2\right)\\
&\le \,\oldgianni{C_9}\left(\|B^\sigma\oldgianni\zk(t)\|\,+\,\|\oldgianni\zk(t)\|\,+\,\|\yk(t)-\by(t)\|_{\Lx4}^2\right),
\end{align*}
\oldgianni{and it follows} the estimate
\oldgianni{%
\begin{align*}
\|\mean(\rk)\|_{L^2(0,T)}\,
\le\,C_{10}\left( 
  \norma\zk_{\L2H}
  + \norma{B^\sigma\zk}_{\L2H}
  + \norma{\yk-\by}_{\L\infty{\Lx4}}^2
\right)\,.
\end{align*} 
}%
In view of \eqref{susi6} \pier{and \eqref{susi4}}, 
\oldgianni{\juerg{and by} recalling \eqref{defnormaAr} and Remark~\ref{Meanvalue}},
this yields that
\oldgianni{%
\Beq
  \norma\rk_{\L2{\VA r}}^2
  \leq C_{11} \left(
    \norma{A^r\rk}_{\L2H}^2
    + \norma{\mean(\rk)}_{L^2(0,T)}^2
  \right) 
  \leq C_{12} \, \norma k_{\L2H}^4 \,.
  \non
\Eeq
}%
In conclusion, the condition \eqref{Frechet} holds true with the choice $\,Z(\lambda)=(C_7+C_{1\oldgianni2})\,\lambda^4$.
With this, the assertion is completely proved.
\Edim

Using the above differentiability result and the fact that $\uad$ is a closed and
convex subset of ${\cal X}$, we can infer from the chain rule via a standard 
argument (which can be omitted here) the following first-order necessary optimality
condition:

\Bcor
\label{FirstNC}
Let the assumptions of \oldgianni{Theorem~\ref{Differentiability}} be satisfied, and assume that $\bar u\in\uad$
with $(\bar\mu,\bar y)={\cal S}(\bar u)$ is a solution to the optimal control problem
{\bf (CP)}. Then it holds the variational inequality
\revis{%
\begin{align}
\label{vug1}
&\alpha_1\iO(\bar y(T)-y_\Omega)\,\xi(T)\,+\,\alpha_2\int_0^T\!\!\iO(\bar y-y_Q)\,\xi
\nonumber\\
& \quad {}
+\,\revis{\alpha_3}\int_0^T\!\!\iO \bar u\,(v-\bar u)\,\ge\,0 \qquad\forall\,v\in\uad,
\end{align}
}%
where $(\eta,\xi)$ is the unique solution to the system \eqref{wlin1}--\eqref{wlin3}
associated with $\,k=v-\bar u$.
\Ecor 

\section{Existence and first-order optimality conditions}
\label{OPTIMUM}
\setcounter{equation}{0}
In this section, we state and prove the main results of this paper. We begin with an existence result.

\Bthm
\label{ExistenceOC}
Suppose that the conditions {\bf (A1)}--{\bf \revis{(A9)}} and {\bf (GB)} are fulfilled. 
Then the optimal control problem {\bf (CP)} has a solution.
\Ethm 
 
\Bdim
We use the direct method. 
To this end, let $\,\{u_n\}\subset \uad\,$ be a minimizing sequence, and let $\,(\mu_n,y_n)={\cal S}(u_n)$, for $n\in\enne$. 
Then the global bounds \eqref{bounds1} and \eqref{bounds2} apply, and there are  some $\,\bu \in\uad$, a pair $\,(\bm,\by)\,$, and some $z\in L^\infty(Q)$, such that,
at least for a subsequence which is again indexed by $n\in\enne$,
\begin{align}
\label{conu}
&u_n\to \bu\quad\mbox{weakly star in \,${\cal X}$,}\\[1mm]
\label{conmy}
&\mu_n\to \bm \quad\mbox{weakly star in \,$L^\infty(0,T;V_A^{2r})$,}\\[1mm]
\label{cony}
&y_n\to \by\quad\mbox{weakly star in \,$W^{1,\infty}(0,T;H)\cap H^1(0,T;V_B^\sigma)$}\,,\\[1mm]
\label{conbeta}
&f_1 '(y_n)\to z\quad\mbox{weakly star in \,$L^\infty(Q)$}. 
\end{align} 
We also observe that standard compactness results (see, e.g. \cite[Sect.~8, Cor.~4]{Simon}) imply that we may without loss of generality  assume that
\begin{equation}
\label{strongcony}
y_n\to y \quad\mbox{strongly in  $\, C^0([0,T];H)\,$ and pointwise a.e. in \,$Q$},
\end{equation}
which yields that $\bar y(0)=y_0$, in particular.
In addition, by {\bf (A4)}, $f_2'$ is Lipschitz continuous on $\erre$, 
which implies that $f_2'(y_n)\to f_2'(\by)$ \,strongly in \pier{$C^0([0,T];H)$}; 
moreover, 
\pier{by the convexity of $f_1$ it turns out that $f_1'$ induces} a maximal monotone graph. It then follows from standard results on maximal monotone operators
(see, e.g., \cite[Prop.~2.2, p.~38]{Barbu}) that $z=f_1'(\by)$. In summary, we have 
\pier{that} $f'(y_n)\to f'(\by)$ weakly \pier{star in $L^\infty(0,T;H)$}. 

\oldgianni{Now, we consider the (equivalent) integrated version of \eqref{weak1}--\eqref{weak3},  
written for $u=u_n$, $y=y_n$, $\mu=\mu_n$, $n\in\enne$, and with time-dependent test functions,
and we pass to the limit as $n\to\infty$}. 
We then obtain the \oldgianni{\juerg{analogous}} formulation for $u=\bar u$, $\mu=\bm$, $y=\by$, 
that is, we have $(\bm,\by)={\cal S}(\bu)$.
But this means that the pair $((\bm,\by),\bu)$ is admissible for the minimization problem {\bf (CP)}. 
By the semicontinuity properties
of the cost functional, it is a minimizer.
\Edim  

\pier{Next, we} aim to establish meaningful first-order necessary optimality conditions by eliminating
the quantities $\eta$ and $\xi$ from \eqref{vug1} by means of the adjoint state variables. 
To this end, we consider the adjoint state system which formally reads
\begin{align}
A^{2r}p-q\,
&= \revis 0
\quad\mbox{in $\,Q$},
\label{adj1}
\\[0.5mm]
-(\dt p+\tau\,\dt q)+B^{2\sigma}q+f''(\by)\,q\,
&=\,\alpha_2(\by-y_Q)
\quad\mbox{in $\,Q$,}
\label{adj2}
\\[0.5mm]
p(T)+\tau\,q(T)\,
&=\,\alpha_1(\by(T)-y_\Omega)
\quad\mbox{in $\,\Omega$.}
\label{adj3}
\end{align}
\revis{Precisely,} 
we consider a variational formulation of the above formal problem.
We recall the definition~\eqref{defVBneg} of $\VB{-\sigma}$
and the embedding $H\subset\VB{-\sigma}$ (see~\eqref{embVBneg})\pier{; let us}
use the simpler notation $\<\cpto,\cpto>$ without indices for the 
duality pairing between $\VB{-\sigma}$ and~$\VB\sigma$.
For the adjoint state~$(p,q)$, we require the following regularity conditions:
\begin{align}
  & p \in \L2{\VA{2r}},
  \label{regp}
  \\
  & q \in \L2{\VB\sigma},
  \label{regq}
  \\
  & p+\tau q \in \H1{\VB{-\sigma}}.
  \label{regdtpq}
\end{align}
\Accorpa\Regadjoint regp regdtpq
The adjoint problem we consider \juerg{then reads as follows:}
\begin{align}
  & (A^r p(t),A^r v) - (q(t),v)
  = 0
  \quad \hbox{\aat\ and every $v\in\VA r$,}
  \label{primaA}
  \\[2mm]
  & - \< \dt(p+\tau q)(t) , v >
  + ( B^\sigma q(t) , B^\sigma v )
  + ( \psi(t) q(t) , v )
  = ( g_2(t) , v )
  \non
  \\
  & \quad \hbox{\aat\ and every $v\in\VB\sigma$,}
  \label{secondaA}
  \\[2mm]
  & (p+\tau q)(T) = g_1, 
  \label{cauchyA}
\end{align}
\Accorpa\Adjoint primaA cauchyA
where, for brevity, we have set 
\Beq
  \psi := f''(\by) , \quad
  g_1 := \alpha_1(\by(T)-y_\Omega)
  \aand
  g_2 := \alpha_2(\by-y_Q) .
  \label{notations}
\Eeq
We have written for convenience the weak form \eqref{primaA},
which still makes sense under the weaker regularity requirement $p\in\L2{\VA r}$.
However, it is \juerg{immediately} seen that such a regularity and \eqref{primaA} imply \eqref{regp} and
\Beq
  q = A^{2r} p,
  \label{strongprimaA}
\Eeq
\revis{i.e., the equation \eqref{adj1}}.

\revis{Solving problem} \Adjoint\ \juerg{requires} some preliminary work.
It is understood that the assumptions {\bf (A1)--\revis{(A9)}} and {\bf (GB)} are in force.
In particular, \juerg{we have that} $\psi\in\LQ\infty$, $g_1\in\Ldue$, and $g_2\in\LQ2$.
First of all, we give an equivalent formulation.

\Bprop
\label{NewAdjoint}
The regularity conditions \Regadjoint\ and problem \Adjoint\ are equivalent 
to \accorpa{regp}{regq}, \eqref{primaA}, and
\Bsist
  && \ioT \bigl( (p+\tau q)(t) , \dt v(t) \bigr) \, dt
  \non
  \\
  && = - \ioT \bigl( B^\sigma q(t) , B^\sigma v(t) \bigr) \, dt
  + \ioT \bigr( g_2(t) - \psi(t) q(t) , v(t) \bigr) \, dt
  + \bigl( g_1 , v(T) \bigr)
  \non
  \\
  && \quad \hbox{for every $v\in\H1H\cap\L2{\VB\sigma}$ satisfying $v(0)=0$}.
  \label{newsecondaA}
\Esist
\Eprop

\Bdim
Before starting, we observe that, for $p\in\L2H$ and $q\in\L2{\VB\sigma}$ with $p+\tau q\in\H1{\VB{-\sigma}}$,
the variational equation \eqref{secondaA} is equivalent to the following integrated version:
\Bsist
  && - \ioT \< \dt(p+\tau q)(t) , v(t) > \, dt 
  \non
  \\
  && = - \ioT \bigl( B^\sigma q(t) , B^\sigma v(t) \bigr) \, dt
  + \ioT \bigr( g_2(t) - \psi(t) q(t) , v(t) \bigr) \, dt 
  \non
  \\
  && \quad \hbox{for every $v\in\L2{\VB\sigma}$}.
  \label{intsecondaA}
\Esist
We also recall an \juerg{integration-by-parts} formula (see, e.g., \cite[Lemma~4.5]{CGSSIAM18}):
if $(\calV,\calH,\calVp)$ is a Hilbert triplet and 
\Beq
  w \in \H1\calH \cap \L2\calV
  \aand
  z \in \H1{\calVp} \cap \L2\calH ,
  \non
\Eeq
then the function $t\mapsto(w(t),z(t))_{\calH}$ is absolutely continuous, and for every $t,t'\in[0,T]$ we have that
\Beq
  \int_{t'}^t \bigl\{
    \bigl( \dt w(s),z(s) \bigr)_{\calH}
    + {}_{\calVp}\< \dt z(s),w(s)>_{\calV}
  \bigr\} \, ds
  = \bigl( w(t),z(t) \bigr)_{\calH}
  - \bigl( w(t'),z(t') \bigr)_{\calH} \,.
  \label{leibniz}
\Eeq
Now, we prove the statement.
We first assume that \Regadjoint\ and \Adjoint\ are valid.
Then, we just have to prove that \eqref{newsecondaA} holds true.
We start from \eqref{intsecondaA}, with any $v\in\H1H\cap\L2{\VB\sigma}$.
By applying~\eqref{leibniz},
we immediately obtain \eqref{newsecondaA} on account of \eqref{cauchyA}.

Conversely, assume that $(p,q)$ satisfies \accorpa{regp}{regq}, \eqref{primaA}, and~\eqref{newsecondaA}.
We prove the (apparently) stronger regularity requirement \eqref{regdtpq} 
and the validity of the formulas \eqref{secondaA} and~\eqref{cauchyA}.
To this end, we observe that, because of the meaning of the Hilbert triplet $(\VB\sigma,H,\VB{-\sigma})$,
the conditions \eqref{regdtpq} and \eqref{secondaA} (or~\eqref{intsecondaA})
are equivalent to the following properties:
$i)$~formula \eqref{newsecondaA} holds for every $v\in C^\infty_c(0,T;\VB\sigma)$ 
($C^\infty$~functions with compact support in~$(0,T)$);
$ii)$~the maps that associates to every $v\in C^\infty_c(0,T;\VB\sigma)$ 
the \rhs\ of \eqref{newsecondaA} (i.e.,~the same as in~\eqref{intsecondaA} since $v(T)=0$)
is~continuous with respect to the topology of~$\L2{\VB\sigma}$.
The former follows from our assumption and it is \sfw\ to see that the latter is satisfied since 
$q\in\L2{\VB\sigma}$.
Hence, both \eqref{regdtpq} and \eqref{intsecondaA} are established
(the latter first for every $v\in C^\infty_c(0,T;\VB\sigma)$ by definition,
then for every $v\in\L2{\VB\sigma}$ by continuity).
At this point,
we can \pier{take \eqref{newsecondaA},} with $v\in\H1H\cap\L2{\VB\sigma}$ 
satisfying $v(0)=0$\pier{, and integrate by parts using \eqref{leibniz}}.
By comparing with \eqref{intsecondaA}, we deduce that
$\bigl( (p+\tau q)(T) - g_1 , v(T) \bigr) = 0$
for every $v\in\H1H\cap\L2{\VB\sigma}$ satisfying $v(0)=0$.
By choosing $v(t)=tv_0$ with any $v_0\in\VB\sigma$ we conclude that \eqref{cauchyA} holds true as~well
since $\VB\sigma$ is dense in~$H$.
\Edim

Thus, we are going to solve the new problem given by the previous proposition.
The case \juerg{$\lambda_1>0$} is easier, since the operator $A^{2r}\in\calL(\VA{2r},H)$ has the inverse 
$A^{-2r}:=(A^{2r})^{-1}\in\calL(H,\VA{2r})$,
so that we can use \eqref{strongprimaA} in order to eliminate~$p$.
Hence, we immediately obtain the following lemma:

\Blem
\label{Soloqpos}
Assume $\lambda_1>0$.
Then, a~pair $(p,q)$ satisfying \accorpa{regp}{regq} solves \eqref{primaA} and \eqref{newsecondaA} if and only if
$p = A^{-2r}q$ with $q$ satisfying
\Bsist
  && q \in\L2{\VB\sigma}, 
  \label{regqpos}
  \\
  && \ioT \bigl( A^{-2r}q(t) + \tau q(t) , \dt v(t) \bigr) \, dt
  \non
  \\
  && = - \ioT \bigl( B^\sigma q(t) , B^\sigma v(t) \bigr) \, dt
  + \ioT \bigr( g_2(t) - \psi(t) q(t) , v(t) \bigr) \, dt
  + \bigl( g_1 , v(T) \bigr)
  \non
  \\
  && \quad \hbox{for every $v\in\H1H\cap\L2{\VB\sigma}$ satisfying $v(0)=0$}.
  \label{equazqpos}
\Esist
\Accorpa\Pblqpos regqpos equazqpos
\Elem

On the contrary, the situation is much more complicated in \juerg{the case when} $\lambda_1=0$.
\juerg{To handle this case, we} adapt the ideas of \cite[Sect.~5]{CGSAMO}.
To this end, we \juerg{have to introduce} some new spaces.
We~set
\Beq
  \Hz := \graffe{v\in H:\ \mean(v)=0} 
  \aand
  \VBz\sigma := \VB\sigma \cap \Hz ,
  \label{defVBz}
\Eeq
and notice that $\H1H\cap\L2{\VBz\sigma}=\H1\Hz\cap\L2{\VBz\sigma}$.
Moreover, we observe that the operators $\Az{2r}$ and $\Az{-2r}$ (see \eqref{defAz}) also satisfy \juerg{that}
\Beq
  \Az{2r} : \Vz{2r} \to \Hz
  \aand
  \Az{-2r} = (\Az{2r})^{-1} : \Hz \to \Vz{2r}
  \quad \hbox{are isomorphisms}.
  \label{isoAz}
\Eeq
Finally, for simplicity, in the next statement and in its proof, we often use the same notation $\phi$
for some real function $\phi\in H^1(0,T)$ and the function $\phi\1\in\H1H$.

\Blem
\label{Soloqzero}
Assume that $\lambda_1=0$.
Then a~pair $(p,q)$ satisfying \accorpa{regp}{regq} solves \eqref{primaA} and \eqref{newsecondaA} if and only~if
\Beq
  p = \pO + \Az{-2r} q,
  \label{calcolop}
\Eeq
with $q$ and $\pO$ given as follows:
\Bsist
  && q \in \L2{\VBz\sigma},
  \label{regqzero}
  \\[1mm]
  && \ioT \bigl( \Az{-2r}q(t) + \tau q(t) , \dt v(t) \bigr) \, dt
  \non
  \\
  && = - \ioT \bigl( B^\sigma q(t) , B^\sigma v(t) \bigr) \, dt
  + \ioT \bigr( g_2(t) - \psi(t) q(t) , v(t) \bigr) \, dt
  + \bigl( g_1 - \mean(g_1)\1 , v(T) \bigr)
  \non
  \\[1mm]
  && \quad \hbox{for every $v\in\H1H\cap\L2{\VBz\sigma}$ satisfying $v(0)=0$,}
  \label{equazqzero}
  \\[2mm]
  && \pO(t) 
  \,=\, \mean(g_1)\,
  + \,\frac 1{|\Omega|} \int_t^T \Bigl\{
    \bigl( g_2(s) , \1 \bigr)    
    - \bigl( B^\sigma q(s) , B^\sigma \1 \bigr)
    - \bigl( \psi(s) q(s) , \1 \bigr)
  \Bigr\} \, ds
  \non
  \\
  && \quad \hbox{for every $t\in[0,T]$}.
  \label{calcolopO}
\Esist
\Accorpa\Pblqzero regqzero equazqzero
\Elem

\Bdim
Assume that $(p,q)$ satisfies \accorpa{regp}{regq} and solves \eqref{primaA} and \eqref{newsecondaA}.
By Proposition~\ref{NewAdjoint} we can also use the previous formulation \Adjoint\ of the adjoint problem.
Testing \juerg{\eqref{primaA}} by $v=\1\in\VB\sigma$ yields
\Beq
  \bigl( q(t),\1 \bigr)
  = \bigl( A^r q(t),A^r \1 \bigr)
  = 0
  \quad \hbox{\aat},
  \non
\Eeq
since $\lambda_1=0$.
Thus, $q$ has zero mean value, and \eqref{regqzero} is a consequence of~\eqref{regq}.
\pier{Moreover, in view of \eqref{regdtpq} it turns out that the function 
$$
%  \pO(t) := 
\giapie{t \mapsto  \mean((p+ \tau q)(t)  ),  \quad t\in [0,T]},
$$
belongs to $H^1(0,T)$, \giapie{and in particular it has a continuous representative (termed exactly as it~is). We set 
\Beq
 \pO(t) := \mean((p+ \tau q)(t)  ),  \quad \hbox{for every } \, t\in [0,T],
  \label{meanp}
\Eeq
and it turns out that 
$$
 \pO(t)  = \mean(p(t)) \quad \aat . 
$$
Therefore, by choosing $v=\1$ in~\eqref{secondaA} and using~\eqref{cauchyA}, 
we also deduce that}}
\Bsist
  && - \,|\Omega| \, \pO'(t)
  + \bigl( B^\sigma q(t) , B^\sigma \1 \bigr)
  + \bigl( \psi(t) q(t) , \1 \bigr)
  = \bigl( g_2(t) , \1 \bigr)
  \quad \aat,
  \non
  \\[1mm]
  && \pO(T) = \mean (g_1) .
  \non
\Esist
Hence, \eqref{calcolopO} immediately follows.
Furthermore, since $A^{2r}\1=0$, we can write \eqref{strongprimaA} in the form
\Beq
  q = A^{2r} (p-\pO) = \Az{2r} (p-\pO), 
  \non
\Eeq
and, owing to the zero mean value property of $q$ once more, we conclude that
\Beq
  p - \pO = \Az{-2r} q, 
  \non
\Eeq
that is, \eqref{calcolop} \juerg{holds true}.
Using this, we compute both sides of \eqref{newsecondaA} 
with zero-mean-value test functions, i.e., $v\in\H1H\cap\L2{\VBz\sigma}$\juerg{,} such that $v(0)=0$.
Since $\pO(t)$ is space independent, $\mean (g_1)$ is a constant,
and $\dt v(t)$ and $v(T)$ have zero mean value,
we~have 
\Bsist
  && \ioT \bigl( (p+\tau q)(t),\dt v(t) \bigr) \, dt
  = \ioT \bigl( (\Az{-2r}q+\tau q)(t),\dt v(t) \bigr) \, dt
  + \ioT \bigl( \pO(t),\dt v(t) \bigr) \, dt
  \non
  \\
  && = \ioT \bigl( (\Az{-2r}q+\tau q)(t),\dt v(t) \bigr) \, dt\,, 
  \non
\Esist
as well as
\Bsist
  && - \ioT \bigl( B^\sigma q(t) , B^\sigma v(t) \bigr) \, dt
  + \ioT \bigr( g_2(t) - \psi(t) q(t) , v(t) \bigr) \, dt
  + \bigl( g_1 , v(T) \bigr)
  \non
  \\
  &&\juerg{=}\, - \ioT \!\bigl( B^\sigma q(t) , B^\sigma v(t) \bigr) \, dt
  + \ioT \!\bigr( g_2(t) - \psi(t) q(t) , v(t) \bigr) \, dt
  + \bigl( g_1 - \mean( g_1) \1 , v(T) \bigr) .
  \non
\Esist
Hence, \eqref{newsecondaA} with such test functions becomes~\eqref{equazqzero}.

Conversely, assume that $p$ fulfils \eqref{calcolop} 
with $q$ satisfying \Pblqzero\ \pier{and with} $\pO$ given by~\eqref{calcolopO}.
First of all, observe that~\eqref{meanp} (which is not required a priori)
still holds as a consequence of~\eqref{calcolop}, since $\Az{-2r}q$ has zero mean value.
Moreover, \eqref{regq} is trivially implied by~\eqref{regqzero}.
We now prove \juerg{the validity of} \eqref{newsecondaA}.
\juerg{To this end, take}
 any $v\in\H1H\cap\L2{\pier{\VB\sigma}}$ with $v(0)=0$ and split $v$ as follows:
\Beq
  v = (v - \phi \1) + \phi \1
  \quad \hbox{where} \quad
  \phi := \mean(v) \,.
  \non
\Eeq
Then, $v-\phi\1\in\H1H\cap\L2{\VBz\sigma}$, and $\,(v-\phi\1)(0)=0$.
Hence, \eqref{equazqzero} yields
\Bsist
  && \ioT \bigl( \Az{-2r}q(t) + \tau q(t) , \dt v(t) - \phi'(t) \, \1 \bigr) \, dt
  \non
  \\
  && = - \ioT \bigl( B^\sigma q(t) , B^\sigma v(t) \bigr) \, dt
  + \ioT \bigr( g_2(t) - \psi(t) q(t) , v(t) \bigr) \, dt
  \non
  \\
  && \quad {}
  + \bigl( g_1 - \mean(g_1)\1 , v(T) \bigr)
  \non
  \\
  && \quad {}
  + \ioT \bigl( B^\sigma q(t) , B^\sigma \1 \bigr) \phi(t) \, dt
  - \ioT \bigr( g_2(t) - \psi(t) q(t) , \1 \bigr) \phi(t) \, dt
  \non
  \\
  && \quad {}
  - \bigl( g_1 - \mean(g_1)\1 , \1 \bigr) \phi(T),
  \non
\Esist
and we note that the last term vanishes.
Now, we observe that $\phi\in H^1(0,T)$ and that $\phi(0)=0$ (since $v(0)=0$).
Thus, we multiply \eqref{calcolopO} by \,$|\Omega|\phi'(t)$,
integrate over $(0,T)$ with respect to~$t$, and perform an integration by parts on the \rhs.
We obtain that
\Bsist
  && \ioT \bigl( \pO(t) , \1 \bigr) \, \phi'(t) \, dt
  \non
  \\
  && = \bigl( \mean(g_1) \1 , \1 \bigr) \phi(T) 
  - \ioT \bigl( B^\sigma q(t) , B^\sigma \1 \bigr) \phi(t)
  + \ioT \bigr( g_2(t) - \psi(t) q(t) , \1 \bigr) \phi(t) \, dt \,.
  \non
\Esist
By summing up, we deduce that
\Bsist
  && \ioT \bigl( \Az{-2r}q(t) + \tau q(t) , \dt v(t) - \phi'(t) \, \1 \bigr) \, dt
  + \ioT \bigl( \pO(t) , \1 \bigr) \, \phi'(t) \, dt
  \non
  \\
  && = - \ioT \bigl( B^\sigma q(t) , B^\sigma v(t) \bigr) \, dt
  + \ioT \bigr( g_2(t) - \psi(t) q(t) , v(t) \bigr) \, dt
  + \bigl( g_1 , v(T) \bigr) .
  \non
\Esist
Notice that the \rhs s of this identity and \juerg{of} \eqref{newsecondaA} coincide.
Thus, it suffices to show that the same happens for the \lhs s.
By also accounting for~\eqref{calcolop}, and noting that 
the mean values of both $\dt v-\phi'\1$ and $p-\pO+\tau q$ vanish, we~have
\Bsist
  && \ioT \bigl( \Az{-2r}q(t) + \tau q(t) , \dt v(t) - \phi'(t) \, \1 \bigr) \, dt
  + \ioT \bigl( \pO(t) , \1 \bigr) \, \phi'(t) \, dt
  \non
  \\
  \separa
  && = \ioT \bigl( p(t) + \tau q(t) - \pO(t) , \dt v(t) - \phi'(t) \, \1 \bigr) \, dt
  + \ioT \bigl( \pO(t) , \1 \bigr) \, \phi'(t) \, dt
  \non
  \\
  \separa
  && = \ioT \bigl( p(t) + \tau q(t) , \dt v(t) \bigr) \, dt
  - \ioT \bigl( p(t) + \tau q(t) , \1 \bigr) \phi'(t) \, dt
  + \ioT \bigl( \pO(t) , \1 \bigr) \, \phi'(t) \, dt
  \non
  \\
  \separa
  && = \ioT \bigl( p(t) + \tau q(t) , \dt v(t) \bigr) \, dt
  - \ioT \bigl( p(t) - \pO(t) + \tau q(t) , \1 \bigr) \phi'(t) \, dt
  \non
  \\
  \separa
  && = \ioT \bigl( p(t) + \tau q(t) , \dt v(t) \bigr) \, dt \,.
  \non
\Esist
This completes the proof.
\Edim

\Blem
\label{Density}
The space $\VBz\sigma$ is dense in~$\Hz$.
In particular, the Hilbert triplet  
\Beq
  (\VBz\sigma,\Hz,\VBz{-\sigma}),
  \quad \hbox{where} \quad
  \VBz{-\sigma} := (\VBz\sigma)^*,
  \non
\Eeq
is meaningful.
\Elem

\Bdim
We assume that $z\in\Hz$ satisfies $(z,v)=0$ for every $v\in\VBz\sigma$
and deduce that $z=0$.
Take any $v\in\VB\sigma$.
Then $v-\mean(v)\mathbf{1}\in\VBz\sigma$, whence $(z,v-\mean(v)\mathbf{1})=0$.
On the other hand, $(z,\mean(v)\mathbf{1})=0$ since $\mean(z)=0$.
Therefore, $(z,v)=0$.
Since this holds for every $v\in\VB\sigma$ and $\VB\sigma$ is dense in~$H$, we conclude that $z=0$.
\Edim

\Blem
\label{Astratto}
Let $(\calV,\calH,\calVp)$ be a Hilbert triplet 
and let $(\cpto,\cpto)$ and $\<\cpto,\cpto>$ be 
the inner product of $\calH$ and the duality pairing between $\calVp$ and~$\calV$, respectively.
Moreover, let $\calA$ and $\calB$ \juerg{satisfy, with suitable positive constants $M$, $\lambda$, and $\alpha$, the following conditions:}
\begin{align}
  & \calA \in \calL(\calH;\calH) 
  \quad \hbox{is symmetric}\, \pier{;}
  \label{hpA}
  \\
  & \calB(t) \in \calL(\calV;\calVp)
  \quad \aat\pier{;}
  \label{hpB}
  \\
  & \hbox{for every $v,w\in\calV$, the function $t\mapsto\<\calA(t)v,w>$ is measurable on $(0,T)$}\pier{;}  \label{pier1}  
  \\
  & \< \calA v,v > \geq \alpha \, \norma v_{\calH}^2
  \quad \hbox{for every $v\in\calH$}\, \pier{;}
  \label{coercA}
  \\
  & \norma{\calB(t)v}_{\calVp} \leq M \, \norma v_{\calV}
  \quad \hbox{\aat\ and every $v\in\calV$}\, \pier{;}
  \label{bddB}
  \\
  & \< \calB(t) v , v > + \lambda \, \norma v_{\calH}^2
  \geq \alpha \, \norma v_{\calV}^2
  \quad \hbox{\aat\ and every $v\in\calV$}\, .
    \label{wcoercB} 
\end{align}
\juerg{Then}, for every $F\in\L2{\calVp}$ and every $\gamma\in\calH$, 
there exists a unique $q\in\L2\calV$ satisfying
\Bsist
  && \ioT \bigl( \calA q(t) , v'(t) \bigr) \, dt
  + \ioT \< \calB(t) q(t) , v(t) > \, dt
  \non
  \\
  && = \ioT \< F(t) , v(t) >
  + \bigl( \gamma , v(T) \bigr)
  \non
  \\
  && \quad \hbox{for every $v\in\H1\calH\cap\L2\calV$ such that $v(0)=0$}.
  \label{astratto}
\Esist
\Elem

\Bdim
The similar forward problem (presented in a slightly different way,
see~also \cite[Lem.~1.1, p.~44]{Lions} for a similar equivalence) 
is~solved in \cite[Thm.~7.1, p.~70]{Lions} 
under even more general assumptions on the structure
(in~particular, there $\calA$~is also allowed to depend on time)
and equivalent assumptions on the data.
\Edim

\Brem
By arguing as we did for Proposition~\ref{NewAdjoint}, one can easily see that
a function $q\in\L2\calV$ solves \eqref{astratto} if and only if it satisfies
\Beq
  q \in \L2\calV , \quad
  \calA q \in \H1{\calVp} , \quad
  - (\calA q)' + \calB u = F,
  \aand
  (\calA q)(T) = \gamma,
  \non
\Eeq
where the abstract equation holds \aet\ in the sense of~$\calVp$
and the final condition is meaningful since $\calA q\in\C0\calVp$.
\Erem

\revis{At this point, we are ready to state a well-posedness result for the adjoint problem, 
i.e., for \juerg{the} system \Adjoint}.
Namely, we have the following theorem.

\Bthm
\label{WellposedAdjoint}
Suppose that the conditions {\bf (A1)}--{\bf \revis{(A9)}} and {\bf (GB)} are fulfilled.
Moreover, assume that $\pier{\bar u}\in\uad$, and let $(\bar\mu,\bar y)=\calS(\bar u)$ be the corresponding state.
Then the adjoint problem \Adjoint\ has a unique solution $(p,q)$ satisfying \Regadjoint.
\Ethm 
 
\Bdim
Thanks to Proposition~\ref{NewAdjoint} and Lemmas~\ref{Soloqpos} and~\ref{Soloqzero},
it is sufficient to establish well-posedness for the sub-problems that involve just~$q$, i.e.,
\Pblqpos\ and \Pblqzero\ in the cases $\lambda_1>0$ and $\lambda_1=0$, respectively.
However, we can unify these problems by seeing both of them as particular cases of a new one.
To this end, we~set
\Bsist
  & \calH := H , \quad
  \calV := \VB\sigma , \quad
  \calA := A^{-2r} + \tau I,
  \aand
  \gamma := g_1,
  & \quad \hbox{if $\lambda_1>0$,}
  \non
  \\
  & \calH := \Hz , \quad
  \calV := \VBz\sigma , \quad
  \calA := \Az{-2r} + \tau I,
  \aand
  \gamma := g_1 - \mean(g_1)\1,
  & \quad \hbox{if $\lambda_1=0$,} 
  \non
\Esist
where $I$ is the identity map of~$\calH$,
and \juerg{we} define $\calB(t)\in\calL(\calV;\calVp)$~by
\Beq
  \< \calB(t) v,w >
  := \bigl( B^\sigma v,B^\sigma w \bigr)
  + \bigl( \psi(t) v,w \bigr)
  \quad \hbox{\aat\ and every $v,w\in\calV$,}
  \non
\Eeq
in both cases (with different meanings of the notations, e.g.,~$\calV$).
Then, each of the problems we have to solve appears in the form~\eqref{astratto}.
It is \juerg{immediately} seen that the \juerg{assumptions} of Lemma~\ref{Astratto} are fulfilled.
In particular, \eqref{bddB} and \eqref{wcoercB} hold since $\psi$ is bounded.
Hence, the lemma provides a unique solution.
\Edim

We conclude with the first-order necessary condition for optimality
expressed in terms of the adjoint state variables.

\Bthm
\label{GoodNC}
Let the assumptions \revis{{\bf (A1)}--{\bf \revis{(A9)}} and {\bf (GB)}} be satisfied, 
and assume that $\bar u\in\uad$ is a solution \revis{to the optimal control problem~{\bf (CP)}}.
Moreover, let $(\bar\mu,\bar y)={\cal S}(\bar u)$ be the corresponding state,
and let $(p,q)$ be the unique solution to the related adjoint problem.
Then the following variational inequality holds true:
\Beq
  \ioT\!\!\iO (q+\revis{\alpha_3}\bu)(v-\bu)\,\ge\,0 
  \quad\forall\,v\in\uad.
  \label{vug2}
\Eeq  
In particular, if $\revis{\alpha_3}\not=0$, the optimal control $\bu$
is the $\pier{\L2H} $-projection of $-q/\revis{\alpha_3}$ on~$\uad$. 
\Ethm

\Bdim
Fix any $v\in\uad$, set $k:=v-\bu$, and consider the solutions $(\eta,\xi)$ and $(p,q)$
to the corresponding linearized system \accorpa{wlin1}{wlin3}
and the adjoint system \Adjoint, respectively.
We test \eqref{wlin1} and \eqref{wlin2} by $p(t)$ and~$q(t)$, respectively.
Then, we add the resulting equalities to each other and integrate over~$(0,T)$.
By recalling the notations~\eqref{notations}, we obtain~that
\Bsist
  && \ioT \bigl\{
    \bigl( \dt\xi(t),p(t) \bigr) 
    + \bigl( A^r\eta(t),A^r p(t) \bigr) 
  \bigl\} \, dt
  \non
  \\
  && \quad {}
  + \ioT \bigl\{
    \bigl( \tau\,\dt \xi(t),q(t) \bigr)
    + \bigl( B^\sigma \xi(t),B^\sigma q(t) \bigr)
    + \bigl( \psi(t)\xi(t),q(t) \bigr)
  \bigr\} \, dt
  \non
  \\
  && = \ioT \bigl( \eta(t)+k(t),q(t) \bigr) \, dt \,.
  \non
\Esist
At the same time, by testing \eqref{primaA} and \eqref{secondaA} by $-\eta(t)$ and~$-\xi(t)$,
summing up and integrating with respect to~$t$, we have~that
\Bsist
  && \ioT \bigl\{
    - \bigl( A^r p(t),A^r \eta(t) \bigr)
    + \bigl( q(t),\eta(t) \bigr)
  \bigr\} \, dt
  \non
  \\
  && \quad {}
  + \ioT \bigl\{
    \< \dt (p+\tau q)(t),\xi(t) >
    - \bigl( B^\sigma q(t),B^\sigma\xi(t) \bigr)
    - \bigl( \psi(t)q(t),\xi(t) \bigr) 
  \bigr\} \, dt
  \non
  \\
  && = - \ioT \bigl( g_2(t),\xi(t) \bigr) \, dt .
  \non
\Esist
At this point, we add \juerg{these equations}  
and notice that several cancellations occur.
\juerg{We are left with the following identity}:
\Bsist
  && \ioT \bigl\{
    \bigl( \dt\xi(t),(p+\tau q)(t) \bigr)
    + \< \dt(p+\tau q)(t),\xi(t) >
  \bigr\} \, dt
  \non
  \\
  && = \ioT \bigl( k(t),q(t) \bigr) \, dt
  - \ioT \bigl( g_2(t),\xi(t) \bigr) \, dt \,.
\Esist
By applying the \juerg{integration-by-parts} formula \eqref{leibniz} to the \lhs,
\juerg{invoking} the Cauchy conditions \eqref{wlin3} and~\eqref{cauchyA},
and rearranging \juerg{terms},
we deduce~that
\Beq
  \bigl( g_1 , \xi(T) \bigr)
  + \ioT \bigl( g_2(t),\xi(t) \bigr) \, dt 
  = \ioT \bigl( q(t),k(t) \bigr) \, dt \,.
  \label{good}
\Eeq
\revis{On the other hand}, 
the inequality \eqref{vug1} given by Corollary~\ref{FirstNC} reads
\Beq
  \bigl( g_1 , \xi(T) \bigr)
  + \ioT \bigl( g_2(t),\xi(t) \bigr) \, dt 
  + \revis{\alpha_3} \ioT \bigl( \bar u(t),k(t) \bigr) \, dt 
  \geq 0 \,.
  \non
\Eeq
By replacing the sum of the first two integrals by the \rhs\ of~\eqref{good}, 
we obtain~\eqref{vug2} and the proof is complete. \pier{Indeed, the last sentence is 
just a consequence of the Hilbert projection theorem, since $\uad $ is a convex and closed 
subset of $\L2H$.} 
\Edim

%%%%%%%%%%%%%%%%%%%%%%%%%%%%%%%%%%%%%%%%%%%%%%%%%%%%%%%%%%%%%%%%%%%%%%%%

\section*{Acknowledgments}
\pier{This research was supported by the Italian Ministry of Education, 
University and Research~(MIUR): Dipartimenti di Eccellenza Program (2018--2022) 
-- Dept.~of Mathematics ``F.~Casorati'', University of Pavia. 
In addition, PC and CG gratefully acknowledge some other 
financial support from the MIUR-PRIN Grant 2015PA5MP7 ``Calculus of Variations'',}
the GNAMPA (Gruppo Nazionale per l'Analisi Matematica, 
la Probabilit\`a e le loro Applicazioni) of INdAM (Isti\-tuto 
Nazionale di Alta Matematica) and the IMATI -- C.N.R. Pavia.

%%%%%%%%%%%%%%%%%%%%%%%%%%%%%%%%%
%% bibliography
%%%%%%%%%%%%%%%%%%%%%%%%%%%%%%%%%

\vspace{3truemm}

\Begin{thebibliography}{10}

\revis{%
\bibitem{ABG}
H. Abels, S. Bosia, M. Grasselli, Cahn--Hilliard equation with nonlocal singular free energies. Ann. Mat. Pura Appl. (4) 
{\bf 194} (2015), 1071-1106.}

\bibitem{AM}
M. Ainsworth, Z. Mao, Analysis and approximation of a fractional
Cahn--Hilliard equation. SIAM J. Numer. Anal. {\bf 55} (2017),
1689-1718.

\bibitem{AkSeSchi}
G. Akagi, G. Schimperna, A. Segatti,
Fractional Cahn--Hilliard, Allen--Cahn and porous medium equations.
J. Differential Equations {\bf 261} (2016), 2935-2985.

\bibitem{HO1}
H. Antil, E. Ot\'arola,
A FEM for an optimal control problem of fractional
powers  of  elliptic  operators. SIAM J. Control Optim. {\bf 53} (2015), 3432-3456.

\bibitem{HO2}
H. Antil, E. Ot\'arola, An a posteriori error analysis 
for an optimal control problem involving the fractional Laplacian.
IMA J. Numer. Anal. {\bf 38} (2018), 198-226.

\bibitem{HAS}
H. Antil, E. Ot\'arola, A. J. Salgado, A space-time fractional optimal control problem:
analysis and discretization.
 SIAM J. Control Optim. {\bf  54} (2016), 1295-1328. 

\bibitem{Barbu}
V. Barbu,
``Nonlinear Differential Equations of Monotone Type in Banach Spaces''.
Springer,
London, New York, 2010.

\bibitem{BoFiVa}
M. Bonforte, A. Figalli, J. L. V\'azquez,
Sharp boundary behaviour of solutions to semilinear nonlocal elliptic equations.
Preprint arXiv:1710.02731 [math.AP] (2017) 1-35.

\bibitem{Brezis}
H. Brezis,
``Op\'erateurs maximaux monotones et semi-groupes de contractions
dans les espaces de Hilbert''.
North-Holland Math. Stud.
{\bf 5},
North-Holland,
Amsterdam,
1973.

\bibitem{CFGS1}
P. Colli, M. H. Farshbaf-Shaker, G. Gilardi, J. Sprekels,
Optimal boundary control of a viscous Cahn--Hilliard system with dynamic boundary condition 
and double obstacle potentials. SIAM J.
Control Optim. {\bf 53} (2015), 2696-2721.

\bibitem{CFGS2}
P. Colli, M. H. Farshbaf-Shaker, G. Gilardi, J. Sprekels,
Second-order analysis of a boundary
control problem for the viscous Cahn--Hilliard equation with dynamic boundary conditions.
Ann. Acad. Rom. Sci. Ser. Math. Appl. {\bf 7}
(2015), 41-66.

\bibitem{CFGS3}
P. Colli, M. H. Farshbaf-Shaker, J. Sprekels,
A deep quench approach to the optimal control of
an Allen--Cahn equation with dynamic boundary conditions and double obstacles.
Appl. Math. Optim. {\bf 71}
(2015), 1-24.

\bibitem{CGRS}
P. Colli, G. Gilardi, E. Rocca, J. Sprekels,
Optimal distributed control of a diffuse interface model
of tumor growth.
Nonlinearity {\bf 30}
(2017), 2518-2546.

\bibitem{CGSAdvan}
P. Colli, G. Gilardi, J. Sprekels,
A boundary control problem for the pure Cahn--Hilliard equation
with dynamic boundary conditions.
Adv. Nonlinear Anal. {\bf 4}
(2015), 311-325.

\bibitem{CGSAIMS}
P. Colli, G. Gilardi, J. Sprekels, 
Distributed optimal control of a nonstandard nonlocal phase field
system. AIMS Mathematics {\bf 1}
(2016), 246-281.

\bibitem{CGSAMO}
P. Colli, G. Gilardi, J. Sprekels,
A boundary control problem for the viscous Cahn--Hilliard equation
with dynamic boundary conditions.
Appl. Math. Optim. {\bf 72}
(2016), 195-225.

\bibitem{CGSEECT}
P. Colli, G. Gilardi, J. Sprekels,
Distributed optimal control of a nonstandard nonlocal phase field
system with double obstacle potential.
Evol. Equ. Control Theory {\bf 6}
(2017), 35-58.

\bibitem{CGSAnnali}
P. Colli, G. Gilardi, J. Sprekels, On a Cahn--Hilliard system with convection
and dynamic boundary conditions. \juerg{Ann. Mat. Pura Appl. (4) (2018), doi:10.1007/s10231-018-0732-1.}

\bibitem{CGSSIAM18}
P. Colli, G. Gilardi, J. Sprekels,
Optimal velocity control of a viscous Cahn--Hilliard system with convection and dynamic boundary conditions.
\juerg{SIAM J. Control Optim. {\bf 56} (2018), 1665-1691.}

\bibitem{CGSConvex}
P. Colli, G. Gilardi, J. Sprekels,
Optimal velocity control of a convective Cahn--Hilliard system with double
obstacles and dynamic boundary conditions: a `deep quench' approach.
 \juerg{To appear in J. Convex Anal. {\bf 26} (2019) (see also Preprint arXiv: 1709.03892
[math. AP] (2017), pp. 1-30.)}

\bibitem{CGS18}
P. Colli, G. Gilardi, J. Sprekels, Well-posedness and regularity for a generalized
fractional Cahn--Hilliard system. Preprint \juerg{arXiv:1804.11290 [math. AP] (2018), pp. 1-35.}

\bibitem{CS}
P. Colli, J. Sprekels, Optimal boundary control of a nonstandard Cahn--Hilliard system 
with dynamic boundary condition and double obstacle inclusions. In:
``Solvability, Regularity,
Optimal Control of Boundary Value Problems for PDEs''  (P. Colli, A. Favini, E. Rocca, 
G. Schimperna, J. Sprekels, eds.), Springer INdAM Series Vol. {\bf 22}, 151-182.

\bibitem{Duan}
N. Duan, X. Zhao,
Optimal control for the multi-dimensional viscous Cahn--Hilliard equation. 
Electron. J. Differential Equations 2015, Paper No. 165, 13 pp.

%\bibitem{GiMiSchi} 
%G. Gilardi, A. Miranville, and G. Schimperna,
%On the Cahn--Hilliard equation with irregular potentials and dynamic boundary conditions.
%{\it Commun. Pure Appl. Anal.\/} 
%{\bf 8} (2009), 881-912.

\bibitem{FGG}
S. Frigeri, C. G. Gal, M. Grasselli, 
On nonlocal Cahn--Hilliard--Navier--Stokes systems in two dimensions. 
J. Nonlinear Sci. {\bf 26}
(2016), 847-893.

\bibitem{FGGS}
S. Frigeri, C. G. Gal, M. Grasselli, J. Sprekels, 
Two-dimensional nonlocal Cahn--Hilliard--Navier--Stokes
systems with variable viscosity, degenerate mobility and singular potential.
WIAS Preprint Series No. 2309, Berlin 2016, \pier{pp.~1-56.}

\bibitem{FGS}
S. Frigeri, M. Grasselli, J. Sprekels, Optimal distributed control of two-dimensional 
nonlocal Cahn--Hilliard--Navier--Stokes systems with degenerate mobility and singular potential.
WIAS Preprint Series No. 2473, Berlin 2018, \pier{pp.~1-30.}

\bibitem{FRS}
S. Frigeri, E. Rocca, J. Sprekels,
Optimal distributed control of a nonlocal Cahn--Hilliard/Navier--Stokes
system in two dimensions. SIAM J. Control Optim.
{\bf 54} (2016), 221-250.

\bibitem{Fukao}
T. Fukao, N. Yamazaki, A boundary control problem for the equation and dynamic
boundary condition of Cahn--Hilliard type. In: ``Solvability, Regularity,
Optimal Control of Boundary Value Problems for PDEs'' (P. Colli, A. Favini, E. Rocca, 
G. Schimperna, J. Sprekels, eds.), Springer INdAM Series Vol. {\bf 22}, 255-280.

\revis{%
\bibitem{GalDCDS}
C.G. Gal, 
On the strong-to-strong interaction case 
for doubly nonlocal Cahn--Hilliard equations.
Discrete Contin. Dyn. Syst.  {\bf 37} (2017),  131-167.
\bibitem{GalEJAM}
C.G. Gal, 
Non-local Cahn--Hilliard equations with fractional dynamic boundary conditions. European J. Appl. Math. {\bf 28} (2017), 736-788. 
\bibitem{GalAIHP}
C.G. Gal, 
Doubly nonlocal Cahn--Hilliard equations.
Ann. Inst. H. Poincar\'e Anal. Non Lin\'eaire {\bf 35} (2018), 357-392.}%

\bibitem{GV}
C. Geldhauser, E. Valdinoci,
Optimizing the fractional power in a model with stochastic PDE constraints.
Preprint 	arXiv:1703.09329v1 [math.AP] (2017), pp. 1-18.

\bibitem{HHCK}
M. Hinterm\"uller, M. Hinze, C. Kahle, T. Keil,
A goal-oriented dual-weighted adaptive finite element approach for the optimal control of a nonsmooth Cahn--Hilliard--Navier--Stokes system.
WIAS Preprint Series No. 2311, Berlin 2016, \pier{pp.~1-27.}

\bibitem{HKW}
M.  Hinterm\"uller,  T.  Keil,  D.  Wegner,
Optimal  control  of  a  semidiscrete  Cahn--Hilliard--Navier--Stokes
system with non-matched fluid densities. SIAM J. Control Optim. {\bf 55} (2018),
1954-1989.

\bibitem{HW1}
M. Hinterm\"uller, D. Wegner,
Distributed optimal control of the Cahn--Hilliard system including the case
of a double-obstacle homogeneous free energy density. SIAM J. Control Optim.
{\bf 50} (2012), 388-418.

\bibitem{HW2}
M. Hinterm\"uller, D. Wegner,
Optimal control of a semidiscrete Cahn--Hilliard--Navier--Stokes system.
SIAM J. Control Optim. {\bf 52} (2014), 747-772.

\bibitem{HW3}
M. Hinterm\"uller, D. Wegner,
Distributed and boundary control problems for the semidiscrete Cahn--Hilliard/Navier--Stokes 
system with nonsmooth Ginzburg--Landau energies. Topological Optimization
and Optimal Transport, Radon Series on Computational and Applied Mathematics
{\bf 17} (2017), 40-63.

\bibitem{Jerome}
J.W. Jerome, ``Approximation of Nonlinear Evolution Systems''. 
Math. Sci. Engrg.  {\bf 164}, Academic Press, Orlando,
1983.

\bibitem{Lions}
J.L. Lions,
\'Equations Diff\'erentielles Op\'erationnelles et Probl\`emes aux Limites.
Grundlehren, Vol.~111, Springer-Verlag, Berlin, 1961.

\bibitem{Medjo}
\oldgianni{T.T.} Medjo, 
Optimal control of a Cahn--Hilliard--Navier--Stokes model 
with state constraints.
J. Convex Anal.
{\bf 22} (2015), 1135-1172.

\bibitem{MiZe}
A. Miranville, S. Zelik, 
Robust exponential attractors for Cahn--Hilliard type equations with singular potentials. 
Math. Methods Appl. Sci. {\bf 27} (2004), 545-582.

\bibitem{PV}
B. Pellacci, G. Verzini,
Best dispersal strategies in spatially heterogeneous environments: optimization of the 
principal eigenvalue for indefinite fractional Neumann problems.
J. Math. Biol. \pier{{\bf 76} (2018), 1357-1386.}

\bibitem{RS}
E. Rocca, J. Sprekels,
Optimal distributed control of a nonlocal convective Cahn--Hilliard equation by
the velocity in three dimensions. SIAM J. Control Optim. {\bf 53} (2015), 1654-1680.

\bibitem{Simon}
J. Simon,
Compact sets in the space $L^p(0,T; B)$.
 Ann. Mat. Pura Appl.~(4) 
{\bf 146}, (1987), 65-96.

\bibitem{SV}
J. Sprekels, E. Valdinoci, A new class of identification problems: optimizing the
fractional order in a nonlocal evolution equation. SIAM J. Control Optim. {\bf 55}
(2017), 70-93.

\bibitem{WN}
Q.-F. Wang, S.-i. Nakagiri: 
Weak solutions of Cahn--Hilliard equations 
having forcing terms and optimal control problems. 
Mathematical models in functional equations (Japanese) (Kyoto, 1999), 
S\={u}rikaisekikenky\={u}sho K\={o}ky\={u}roku No. 1128 (2000), 172--180.

\bibitem{ZL1}
 X. P. Zhao, C. C. Liu,
Optimal control of the convective Cahn--Hilliard equation. Appl. Anal.
{\bf 92} (2013), 1028-1045.

\bibitem{ZL2}
X. P. Zhao, C. C. Liu,
Optimal control of the convective Cahn--Hilliard equation in 2D case. Appl.
Math. Optim. {\bf 70}
(2014), 61-82.

\bibitem{Z}
J. Zheng, Time optimal controls of the Cahn--Hilliard equation with 
internal control. Optimal Control Appl. Methods {\bf 36}
(2015), 566-582.

\bibitem{ZW}
J. Zheng, Y. Wang, Optimal control problem for Cahn--Hilliard equations with 
state constraint. J. Dyn. Control Syst. {\bf 21} (2015), 257-272.

\End{thebibliography}

\End{document}

%%%%%%%%%%%%%%%%%%%%%%%%%%%%%%%%%%%%%%%%%%%%%